\documentclass[aop]{imsart}

\RequirePackage{amsthm,amsmath,amsfonts,amssymb}
\RequirePackage[numbers]{natbib}
\RequirePackage[colorlinks,citecolor=blue,urlcolor=blue]{hyperref}
\RequirePackage{graphicx}

\startlocaldefs
\theoremstyle{plain}

\newtheorem{theorem}{Theorem}[section]
\newtheorem{lemma}[theorem]{Lemma}
\theoremstyle{remark}

\def\be{\color{black}}

\endlocaldefs

\begin{document}

\begin{frontmatter}
\title{Perturbations of Parabolic Equations and Diffusion Processes with Degeneration: Boundary Problems, Metastability, and 
Homogenization}
\runtitle{Perturbations of Diffusion Processes with Degeneration}

\begin{aug}
\author[A]{\fnms{Mark} \snm{Freidlin}\ead[label=e1]{mif@umd.edu}} \and
\author[B]{\fnms{Leonid} \snm{Koralov}\ead[label=e2]{koralov@umd.edu}}
\address[A]{Dept of Mathematics, University of Maryland,
College Park, MD 20742,
\printead{e1}}

\address[B]{Dept of Mathematics, University of Maryland,
College Park, MD 20742,
\printead{e2}}
\end{aug}

\begin{abstract}
We study diffusion processes that are stopped or reflected on the boundary of a domain. The generator
of the process is assumed to contain two parts: the main part that degenerates on the boundary in a direction orthogonal to
the boundary and a small non-degenerate perturbation. The behavior of such processes determines the stabilization of solutions to the corresponding parabolic equations with a small parameter. Metastability effects arise in this case:
the asymptotics of solutions, as the size of the perturbation tends to zero, depends on the time scale. Initial-boundary value
problems with both the Dirichlet and the Neumann boundary conditions are considered. 
We also consider periodic homogenization for operators with degeneration. 
\end{abstract}



\end{frontmatter}

\section{Introduction and the  main results} \label{intro}

We will study the asymptotic behavior of solutions to parabolic PDEs with the operator that contains two
parts: the main part that degenerates on the boundary of a domain or (in the section on homogenization) on a collection 
of hypersurfaces and
a small non-degenerate perturbation. PDEs with degeneration on the boundary (as a rule, without a perturbation) were considered in 
\cite{Fich}, \cite{Has}, \cite{OR}, \cite{OR2} \cite{F85}  (see also references therein).  Our results include an analysis of boundary
behavior for degenerate processes and corresponding PDEs. This allows us to describe metastability for perturbed systems. In PDE terms,
this means that we obtain the asymptotics of solutions when time $t(\varepsilon)$ tends to infinity as the perturbation size 
$\varepsilon$ tends to zero.  The value of the limit may depend on how the point $(1/\varepsilon, t(\varepsilon))$ approaches infinity. In probabilistic terms, our results provide a description of the long-time behavior of randomly perturbed degenerate stochastic systems that
is similar to the theory of metastability for randomly perturbed dynamical systems (see Chapter 6 of \cite{FW}).

 Using probabilistic representation for solutions of parabolic PDEs, the problems that we consider  can be formulated and solved in  terms of the asymptotic behavior of the corresponding diffusions. We start with the probabilistic description.

Consider a diffusion process $X^{x}_t$ that satisfies the stochastic differential equation
\[
d X^{x}_t = v_0( X^{x}_t) dt + \sum_{i=1}^d v_i( X^{x}_t) \circ d W^i_t,~~~~~X^{x}_0 = x \in \mathbb{R}^d,
\]
where $W^i_t$ are independent Wiener processes and $v_0,...,v_d$ are $C^3(\mathbb{R}^d)$ vector fields. The Stratonovich form is convenient here since it allows
one to provide a coordinate-independent description of the process. The generator of the process, when applied to functions that
are bounded and continuous together with the first and second partial derivatives, is the operator
\begin{equation} \label{opl}
Lu = L_0 + \frac{1}{2} \sum_{i=1}^d L_i^2,
\end{equation}
where $L_i u = \langle v_i, \nabla u \rangle$ is the operator of differentiation along the vector field $v_i$, $i =0,...,d$. 

Let $D \subset \mathbb{R}^d$ be a bounded connected domain. We assume that $\partial D = S_1
\bigcup ... \bigcup S_m$, where $S_1, ... , S_m$ are $C^4$-smooth $(d-1)$-dimensional non-intersecting surfaces (manifolds without boundaries). We will assume that each of the surfaces is invariant for the process and that the diffusion restricted to a single surface is an ergodic process. We also assume that
the diffusion matrix is non-degenerate inside $D$.   

For each $k = 1,...,m$ and each $x \in S_k$, we define $T(x)$ to be the tangent space to $S_k$ at $x$. We
will assume that: 

(a) ${\rm span} (v_0(x),v_1(x),...,v_d(x)) = {\rm span} (v_1(x),...,v_d(x)) = T(x)$ for $x \in S_1 \bigcup ... \bigcup S_m$;

(b) ${\rm span} (v_1(x),...,v_d(x)) = \mathbb{R}^d$ for $x \in D$. 
\\
This is just a more convenient (and slightly stronger) way to express the assumptions that we already made: (a)~implies that the surfaces are invariant, and the process is ergodic on each surface; 
(b)~means that the diffusion matrix is non-degenerate inside the domain. Under these conditions, the process $X^{x}_t$ starting at $x \in D$
does not reach $\partial D$ in finite time.

 Next, we perturb the process $X^x_t$ by a small non-degenerate diffusion. 
The resulting process $X^{x,\varepsilon}_t$ satisfies, for $t \leq \tau^{x,\varepsilon} = \inf\{s: X^{x,\varepsilon}_s \in \partial D\}$, the following stochastic differential equation:
\begin{equation} \label{peqn1}
d X^{x,\varepsilon}_t = (v_0 + \varepsilon^2   \tilde{v}_0) (X^{x,\varepsilon}_t ) dt + \sum_{i=1}^d v_i(X^{x,\varepsilon}_t ) \circ d W^i_t + \varepsilon 
\sum_{i =1}^d \tilde{v}_i (X^{x,\varepsilon}_t ) \circ d \tilde{W}^i_t,~~~X^{x,\varepsilon}_t = x \in \overline{D},
\end{equation}
where $\tilde{W}^i_t$ are independent Wiener processes (also independent of all $W^i_t$), and $\tilde{v}_0,...,\tilde{v}_d$ are 
$C^3(\mathbb{R}^d)$ vector fields.   In order to make our assumption on the 
non-degeneracy of the
perturbation more precise, we state it as follows: 

(c) ${\rm span} (\tilde{v}_1(x),...,\tilde{v}_d(x)) = \mathbb{R}^d$ for $x \in \overline{D}$. 
\\

The process $X^{x,\varepsilon}_t$ is approximated well by $X^{x}_t$ on finite time intervals as $\varepsilon \downarrow 0$. However, due to the presence of the small non-degenerate component,  
$X^{x,\varepsilon}_t$ can reach the surfaces $S_k$, $1 \leq k \leq m$, that are inaccessible for $X^{x}_t$. For $t > \tau^{x,\varepsilon}$,
 $X^{x,\varepsilon}_t$ is defined as the process with co-normal  reflection on the  boundary (the co-normal is obtained by applying the diffusion matrix to the inward-pointing normal and normalizing the resulting vector - 
see, e.g., Sections 1.6 and 2.5 of \cite{F85}; in our case, the co-normal is independent of $\varepsilon$, as explained below). 
If stopped when it hits the boundary, it results in the process  
$X^{x,\varepsilon}_{t \wedge \tau^{x,\varepsilon}}$.

 The reflected and the stopped processes are suitable for studying the initial-boundary value problems with the Neumann  and the  Dirichlet boundary conditions, respectively.  
Inside the domain, in both cases, the generator of the process is $L^\varepsilon = L + \varepsilon^2  \tilde{L}$, with 
\[
\tilde{L}u = \tilde{L}_0 + \frac{1}{2} \sum_{i=1}^d \tilde{L}_i^2,
\]
where $\tilde{L}_i$ is the operator of differentiation along the vector field $\tilde{v}_i$, $i =0,...,d$.

 In order to describe the behavior of $X^{x,\varepsilon}_t$ at times that grow as $\varepsilon \downarrow 0$, we first need to understand the behavior of the process near the surfaces $S_1,...,S_n$. Each surface will be classified
as either attracting or repelling for $X^{x}_t$, depending (roughly speaking) on whether 
$\mathrm{P} (\lim_{t \rightarrow \infty} {\rm dist}(X^x_t, S_k) = 0) > 0$  for each $x \in D$ (attracting) or 
$\mathrm{P} (\lim_{t \rightarrow \infty} {\rm dist}(X^x_t, S_k) = 0) = 0$  for each $x \notin S_k$ (repelling). 
Each surface $S_k$ carries a unique invariant
probability measure for the unperturbed process. This measure will be denoted by $\pi_k$. If 
the boundary of $D$ is a union of repelling surfaces, then $X^x_t$ has, in addition, an invariant measure $\mu$ on $D$ that is 
absolutely continuous with respect to the Lebesgue measure. Yet another collection of measures will be of importance. Let 
$\nu^{x, \varepsilon}$ be the measure on $\partial D$ induced by $X^{x,\varepsilon}_{\tau^{x,\varepsilon}}$ and $\nu^{x, \varepsilon}_k$ be the 
probability measure on $S_k$ obtained by normalizing the restriction of $\nu^{x, \varepsilon}$ to $S_k$. It can be seen that, 
for each $x \in D$, there is a weak limit $\nu_k = \lim_{\varepsilon \downarrow 0} \nu^{x, \varepsilon}_k$, and the limit
does not depend on $x$. These measures can be introduced in a slightly different but equivalent way: by starting the process
$X^{x,\varepsilon}_t$ in a neighborhood of $S_k$, conditioning it on not leaving a slightly larger neighborhood of $S_k$, and then taking the limit, as $\varepsilon \downarrow 0$, of the measures induced by the stopped process.  Note that the measures $\pi_k$ and $\nu_k$ are, in general, different. To see this, consider the example where $D$ is the unit disc and the boundary has a single component $S_1$ (the unit circle in $\mathbb{R}^2$).
Suppose that $L$ is rotation-invariant. Then $\pi_1$ is the uniform distribution on $S_1$. If, in polar coordinates in a
neighborhood of $S^1$, the perturbation has the form $\tilde{L} = \partial^2_\varphi + g(\varphi) \partial^2_r$ with $g$ that is not constant,
then it is possible to show that the measure $\nu_1$ will be not uniform. Roughly speaking, its density with respect to the uniform measure will be larger near the values of $\varphi$ corresponding to the large values of $g$ (the precise construction of the measures $\nu_k$ is found in Section~\ref{neigh}).

Stated in probabilistic terms, the goal of this paper is to describe the behavior of the  processes $X^{x,\varepsilon}_t$   and 
$X^{x,\varepsilon}_{t \wedge \tau^{x,\varepsilon}}$ 
at 
times $t(\varepsilon)$ that tend to infinity as $\varepsilon \downarrow 0$.
Roughly speaking, one can divide the quadrant $(0,\infty) \times (0, \infty)$ into a finite number
of domains such that $X^{x,\varepsilon}_t$ ($X^{x,\varepsilon}_{t \wedge \tau^{x,\varepsilon}}$)  has a limiting distribution (that depends on the initial point)
when $(1/\varepsilon, t(\varepsilon))$ approaches infinity without leaving a given domain. For different domains,
these limits are different. They are referred to as metastable distributions.  These metastable distributions will be seen to be linear combinations of the measures $\pi_1,...,\pi_m$, $\nu_1,...,\nu_m$, and $\mu$ (or a certain subset of this collection,  depending on whether
we consider the reflected or the stopped process and whether any of the surfaces are attracting). 
\\

Let us briefly discuss a more general approach, proposed in \cite{Finf}, that is 
applicable to various problems concerning the long time behavior of perturbed systems. 
Namely, in order to describe the behavior of $X^{x,\varepsilon}_{t(\varepsilon)}$   ($X^{x,\varepsilon}_{t(\varepsilon) 
\wedge \tau^{x,\varepsilon}}$)  at  various time scales $t(\varepsilon)$, one should consider  the simplex ${\mathcal{M}}$ of invariant probability measures of the Markov family of
unperturbed process $X^x_t$, $x \in \overline{D}$. 
The simplex ${\mathcal{M}}$ is the convex envelope of the set $\mathcal{M}_{\rm erg}$ of ergodic invariant probability measures 
of $X^x_t$. Thus, each element of ${\mathcal{M}}$ can be associated with a probability distribution on $\mathcal{M}_{\rm erg}$. In our case, the set $\mathcal{M}_{\rm erg}$ is as follows: 

(a) $\mathcal{M}_{\rm erg}  = \{\pi_1$,...,$\pi_m\}$ for the case with the reflecting boundary if some of the components are attracting.

(b) $\mathcal{M}_{\rm erg}  = \{\pi_1$,...,$\pi_m, \mu\}$ for the case with the reflecting boundary if $S_1,...,S_m$ are repelling.

(c) $\mathcal{M}_{\rm erg}  = \{\delta_y,y \in \partial D\}$ for the case with the absorbing boundary if some of the components are attracting.
In the case of the absorbing boundary, the notion of the unperturbed process 
needs to be modified so that $X^x_t = x$ for all $t$ if $x \in \partial D$, and thus each point of the boundary supports an invariant probability measure.

(d)  $\mathcal{M}_{\rm erg}  = \{\delta_y,y \in \partial D, \mu\}$ for the case with the absorbing boundary if $S_1,...,S_m$ are repelling.

 Let $\eta^x$ be the limiting distribution of $X^x_t$ as $t \rightarrow \infty$. It is possible to show that
such a limiting distribution exists in our case. 
Consider the mapping $\Pi: \overline{D} \rightarrow {\mathcal{M}}$, where $\Pi(x)  = \eta^x$.  
In a broad class of asymptotic problems, including the one considered here,  the behavior of 
$\Pi (X^{x,\varepsilon}_{t(\varepsilon)})$ (and also of $\Pi(X^{x,\varepsilon}_{t(\varepsilon) 
\wedge \tau^{x,\varepsilon}})$, in our case) can be approximated, for small $\varepsilon$, across a variety of time scales, by an $\varepsilon$-dependent
path  in $\mathcal{M}$. The approximating process has the initial distribution $\eta^x$.  
Note that in cases (c) and (d), the measures $\pi_1,...,\pi_m, \nu_1,...,\nu_m$  belong to the span of $\{\delta_y,y \in \partial D\}$ and thus are the elements of $\mathcal{M}$.
\\

Now we can discuss the PDE interpretation of the results on metastability.
We will consider the first initial-boundary value problem: 
\begin{equation} 
\begin{split} \label{direq2} 
&\frac{\partial u^\varepsilon (t,x)}{\partial t} = L^\varepsilon u^\varepsilon (t, x),~~t > 0, x \in D; \\
u^\varepsilon(0, x)  =  &
g(x),~~x \in D;~~~~u^\varepsilon(t,x) = \psi(x),~~t > 0, x \in \partial D,
\end{split}
\end{equation}
where $g \in C(\overline{D})$, $\psi \in C(\partial D)$, and the second initial-boundary value problem:
\begin{equation}
\begin{split} 
 \label{direq3}
& \frac{\partial  u^\varepsilon  (t,x)}{\partial t} = L^\varepsilon u^\varepsilon (t, x),~~t > 0, x \in D; \\
u^\varepsilon(0, x)  = &
g(x),~~x \in D;~~~~\frac{\partial u^\varepsilon(t,x)}{\partial n(x)} = 0,~~t > 0, x \in \partial D,
\end{split}
\end{equation}
where  $n (x)$  is the co-normal (with respect to $L^\varepsilon$) to $\partial D$ at $x$. Since
$L$ degenerates in the direction orthogonal to the boundary, $n(x)$ coincides with the co-normal with respect to $\tilde{L}$  
and does not depend on $\varepsilon$. In fact, as will be clear from our arguments, our results 
will also hold for
any sufficiently smooth field $n$ that is not tangent to the boundary, i.e., the oblique derivative problem can be treated in the same way
as  the second initial-boundary value problem.

The solutions of initial-boundary value problems (\ref{direq2}) and (\ref{direq3}) can be written as expectations of certain functionals of the
corresponding diffusion processes that depend on the parameter $\varepsilon$, and the dependence of 
the limit of $u^\varepsilon(t(\varepsilon),x)$ on the asymptotics of $t(\varepsilon)$ is a
manifestation of metastability for the underlying diffusion.

With each attracting component $S_k$, we can associate a number $\gamma_k >0$ such
that the time it takes the process $X^{x,\varepsilon}_t$ (with reflection) starting at $x \in S_k$ to exit a fixed neighborhood of $S_k$ is of order $\varepsilon^{-\gamma_k}$ when $\varepsilon \downarrow 0$. 
If $S_k$ is repelling, we can associate a number ${\gamma}_k < 0$ to   it
such that the time it takes the process $X^{x,\varepsilon}_t$, $x \in D$, conditioned on exiting the domain through $S_k$,  to reach $S_k$ is of order $\varepsilon^{ {\gamma}_k}$.  
 The numbers $\gamma_k$ can be
found by solving a certain non-linear spectral problem that involves the operator $L$ restricted to $S_k$ and the leading terms of the coefficients
of the operator near~$S_k$. In order to avoid unnecessary technical details, we will assume that $\gamma_k \neq 0$. Without loss of generality, we can assume that the components of the boundary
$S_1,..., S_m$ are labeled in such an order that $\gamma_1 \geq  ... \geq  \gamma_m$. Let $\overline{m}$ be such 
that $\gamma_{\overline{m}} > 0 > \gamma_{\overline{m}+1}$ (we put $\overline{m}  = 0$ if $\gamma_1 < 0$ and $\overline{m} = m$ if $\gamma_m 
> 0$).  

Let us discuss the asymptotics of the first initial-boundary value problem. We will show that, under appropriate assumptions on 
$t(\varepsilon)$, the limit $\lim_{\varepsilon \downarrow 0} u^\varepsilon(t(\varepsilon), x)$
exists and is equal to a certain linear combination of  $\int_{S_k} g  d \pi_k$, $\int_{S_k} \psi d \nu_k$, $1 \leq k \leq m$, and
$\int_D g d \mu$. The coefficients in this linear combination depend on $x \in D$ and on the time scale (but not on
$g$ and $\psi$), and can be calculated as solutions of the appropriate Dirichlet problems for the operator $L$.
Let us formulate a theorem in which the key ingredients  (the measures  
$\pi_1,...,\pi_m, \nu_1,...,\nu_m, \mu$, the coefficients in front of those measures, and the exponents $\gamma_1,...,\gamma_k$) 
were discussed above but will be defined rigorously  in the subsequent sections. 

\begin{theorem} \label{mnt1} Under the above assumptions on the domain $D$ and the operator $L^\varepsilon$, the asymptotics of
solutions to the initial-boundary value problem (\ref{direq2}) is as follows. 

(a) If at least one component of the boundary is attracting (i.e., $\gamma_1 > 0$, which is equivalent to $\overline{m} \geq 1$), then there exist  $p^x_k \geq 0$,  $1 \leq k \leq \overline{m}$, that are continuous 
functions of $x \in D$,  with $\sum_{k =1}^{\overline{m}}p^x_k = 1$, such that 
  \[
	\lim_{\varepsilon \downarrow 0} u^\varepsilon(t(\varepsilon), x) = \sum_{k = 1}^{\overline{m}} p^x_k \int_{S_k} g  d \pi_k,~~
	if~~~1 \ll t(\varepsilon) \ll |\ln(\varepsilon)|,
	\]
  \[
	\lim_{\varepsilon \downarrow 0} u^\varepsilon(t(\varepsilon), x) = \sum_{k = 1}^{\overline{m}} p^x_k \int_{S_k} \psi  d \nu_k,~~
	if~~  t(\varepsilon) \gg |\ln(\varepsilon)|.
	\]
	
(b) If all the components of the boundary are repelling (i.e., $\gamma_1 < 0$, which is equivalent to $\overline{m} = 0$), then 
\[
	\lim_{\varepsilon \downarrow 0} u^\varepsilon(t(\varepsilon), x) =  \int_{D} g  d \mu,~~
	if~~~1 \ll t(\varepsilon) \ll \varepsilon^{\gamma_1},
\]
\[
	\lim_{\varepsilon \downarrow 0} u^\varepsilon(t(\varepsilon), x) =  \int_{S_1} \psi  d \nu_1,~~
	if~~~ t(\varepsilon) \gg \varepsilon^{\gamma_1}.
\]
\end{theorem}
The condition $t(\varepsilon) \ll |\ln(\varepsilon)|$ appears in the first case because it takes time of order $|\ln(\varepsilon)|$ for the process $X^{x,\varepsilon}_t$ to reach  the boundary  $\partial D$  if there is at least one attracting component. Thus the asymptotic 
behavior $u^\varepsilon(t(\varepsilon), x)$  switches at times of order $|\ln(\varepsilon)|$. If all $S_k$ are repelling,  it takes
time of order $\varepsilon^{\gamma_1}$ for the process to reach the boundary, and the switch happens at the time scale $\varepsilon^{\gamma_1}$.

Now let us formulate a theorem on the second initial-boundary value problem (\ref{direq3}). Here, in some cases,
the switch in the asymptotic behavior of $u^\varepsilon(t(\varepsilon), x) $ happens at several time scales.
\begin{theorem} \label{mnt2}
Under the above assumptions on the domain $D$ and the operator $L^\varepsilon$, the asymptotics of
solutions to the initial-boundary value problem (\ref{direq3}) is as follows. 

(a) If at least one component of the boundary is attracting (i.e., $\gamma_1 > 0$), then there exist 
 $p^{x,l}_k \geq 0$,   $1 \leq l \leq \overline{m}$, $1 \leq k \leq l$,  that are continuous 
functions of $x \in D$, with $\sum_{k =1}^{l}p^{x,l}_k = 1$, such that 
 \[
	\lim_{\varepsilon \downarrow 0} u^\varepsilon(t(\varepsilon), x) = \sum_{k = 1}^{\overline{m}} p^{x,\overline{m}}_k 
	\int_{S_k} g  d \pi_k,~~
	if~~~1 \ll t(\varepsilon) \ll \varepsilon^{-\gamma_{\overline{m}}},
	\]
  \[
	\lim_{\varepsilon \downarrow 0} u^\varepsilon(t(\varepsilon), x) = \sum_{k = 1}^{l} p^{x,l}_k \int_{S_k} g  d \pi_k,~~
	if~~~\varepsilon^{-\gamma_{l+1}} \ll t(\varepsilon) \ll \varepsilon^{-\gamma_{l}}~~with~~1 \leq l < \overline{m},
	\]
  \[
	\lim_{\varepsilon \downarrow 0} u^\varepsilon(t(\varepsilon), x) =  \int_{S_1} g  d \pi_1,~~
	if~~~ t(\varepsilon) \gg \varepsilon^{-\gamma_{1}}.
	\]

(b) If all the components of the boundary are repelling (i.e., $\gamma_1 < 0$), then 
\[
	\lim_{\varepsilon \downarrow 0} u^\varepsilon(t(\varepsilon), x) =  \int_{D} g  d \mu,~~
	if~~~ t(\varepsilon) \gg 1.
\]
\end{theorem}

Observe that the measures $\pi_1,...,\pi_m$ and $\mu$ do not depend on the coefficients of the perturbation $\tilde{L}$. Neither do 
the coefficients $p^{x,l}_k$ in the representation of the solution, as will follow from the arguments below. Thus the asymptotic
behavior of 
solutions to the second initial-boundary value problem (\ref{direq3}) (for the time scales satisfying
the assumptions of Theorem~\ref{mnt2}) is determined exclusively by the unperturbed operator $L$ and does not depend on the 
perturbation $\tilde{L}$. The same can not be said about the solutions to the first initial-boundary value problem (\ref{direq2})
since the measures $\nu_k$ do depend on the coefficients of~$\tilde{L}$. 

The probabilistic analogues of the PDE 
results (which imply Theorems~\ref{mnt1} and \ref{mnt2}) will be stated and proved in Section~\ref{mdis},
and the proofs will be based on the arguments in Sections~\ref{anr}-\ref{atpr}.

Let us mention some of the main steps involved. By retaining only the leading terms for the drift and diffusion coefficients of $L$ near 
$S_k$,  we can get an operator that is  homogeneous in the  direction co-normal to $S_k$. In Section~\ref{anr}, we will state and prove a relatively simple but crucial lemma (Lemma~\ref{spec}) on a non-linear spectral problem that
involves the generator of the process restricted to $S_k$ and the leading terms of $L$ in the co-normal direction. 
In Section~\ref{neigh}, we analyze the behavior of the process $X^{x,\varepsilon}_t$ in the vicinity of $S_k$. Due to near-homogeneity 
of $L$ and a certain averaging effect along $S_k$ related to Lemma~\ref{spec}, the process $X^{x,\varepsilon}_t$  is effectively self-similar
in the direction co-normal to $S_k$ at certain spatial scales.
In particular, the exponents $\gamma_k$ that explain the scaling of the transition times near $S_k$ are provided in Lemma~\ref{spec}. 
In Section~\ref{atpr}, we study the asymptotics, as $\varepsilon \downarrow 0$, of the transition probabilities for the process $X^{x,\varepsilon}_t$ between different components of the boundary. In the case when all
the components are attracting, the asymptotics can be expressed in terms of the process $X^x_t$ conditioned on not reaching one of
the components (see the Remark at the end of Section~\ref{atpr}). In Section~\ref{mdis}, we  combine the results on the local (Section~\ref{neigh})
and the global (Section~\ref{atpr}) behavior of the process to prove the main results, stated in probabilistic terms. In Section~\ref{homog1},
we discuss homogenization for  processes with degeneration on a periodic array of hypersurfaces. 
Due to metastability, different effective (homogenized) processes may appear, depending on the time scale. 
\\

By our assumptions on the vector fields $v_0,...,v_d$, the boundary  $\partial D = S_1
\bigcup ... \bigcup S_m$ is inaccessible in finite time with probability one for the process $X^{x}_t$, $x \in D$. If we change
assumptions on $v_0$ to allow $\partial D$ to be accessible in finite time, then, under mild additional assumptions, the 
solution of the Dirichlet problem for the operator $L$ exists and is unique (\cite{F85}, Chapter 3), and the solution of the first 
initial-boundary value problem for the perturbed equation converges to the unique solution of the non-perturbed problem.
 Thus the case considered in this paper, when both the drift and diffusion degenerate in the direction orthogonal to the boundary, is,
in a sense, the
most interesting.  Still, there remain open questions for both first and second initial-boundary value problems, particularly in the case when
the boundary is inaccessible with $v_0$ forming an acute angle with the inward-pointing normal on the boundary. Such problems will be 
considered in a subsequent paper.

Let us also note a related series of questions concerning random perturbations of dynamical systems where the unperturbed vector
field is tangent to the boundary. In the situation where, for each $x \in D$, the trajectory  starting at $x$ is attracted to a set
inside $D$,  the problem was considered
by M. Day (\cite{D1}, \cite{D2}, \cite{D3}).  It was shown that the  exit measure for the perturbed process can be related to the invariant measure for the process with reflection, which allows to investigate the limit of the exit measures as the size of
the perturbation goes to zero. Similar results, without rigorous justification, were earlier obtained in \cite{MS}. The qualitative behavior of such randomly perturbed dynamical systems is similar, to some extent, to our case of repelling boundary (however, without diffusion in the unperturbed system). In the situation when the unperturbed process has a first integral (and averaging can be used), 
the behavior of the perturbed process and the exit problem were 
considered in \cite{H2}, \cite{FW}. 


Finally, let us mention that a more general problem appears if one forgoes the assumption that $S_k$ serve as boundary components (and thus
are hypersurfaces in $\mathbb{R}^d$). One could instead consider the process $X^{x,\varepsilon}_t$ on a bounded domain $D$ with reflection on $\partial D$ such
that the unperturbed process $X^x_t$ degenerates on a collection of surfaces $S_k \subset D$ of various dimensions (the boundary of $D$ may also contain components of various dimensions). Analyzing the local 
behavior of the process becomes more complicated since the component of $X^{x,\varepsilon}_t$ in the direction(s) transversal to $S_k$
can now be a multi-dimensional process. However, ideas similar to those presented in this paper can still be employed. 

In a recent paper \cite{FK-MC}, we demonstrated metastability for families of parameter-dependent Markov chains, provided that the transition rates satisfy a natural 
assumption (that we called complete asymptotic regularity). Similar results were obtained in
\cite{LX} and  \cite{BL}. In an upcoming paper, we will prove an abstract result on metastability for parameter-dependent Markov renewal processes that will be applicable to various continuous-time systems. Namely, by introducing a sequence of stopping times that correspond to
a continuous-time process 
reaching a vicinity of metastable state, one obtains a Markov renewal process. Examples of such continuous-time processes include the randomly perturbed processes with degeneration 
discussed above, random perturbations of dynamical systems with multiple stable equilibria (\cite{FW}), or motion along heteroclinic
networks  (\cite{Ba}, \cite{BCP}). Of course, application of an abstract result to a given continuous-time system will require, 
each time, verifying
an analogue of the complete asymptotic regularity conditions.

\section{Assumptions on the coefficients, and the structure of the operators $L$ and $L^\varepsilon$ near the boundary} \label{anr}

We assume that the boundary $\partial D$ is $C^4$-smooth, $v_0,...,v_d, \tilde{v}_0,...,\tilde{v}_d \in C^3(\overline{D})$, and
that (a)-(c) hold. Let us specify how the coefficients of the process $X^x_t$ degenerate on $\partial D$. 
Roughly speaking, while there is no diffusion or drift across the boundary, the diffusion should degenerate in a generic way, i.e., the  derivative of the
diffusion coefficient in the direction transversal to $\partial D$ should be non-zero. Let us make this assumption more precise.  

Fix $k \in \{1,...,m\}$. Recall that $ {n}$ is the field of co-normal vectors (with respect to $L^\varepsilon$) on $\partial D$.   
For each $x \in \overline{D}$ 
in a sufficiently small neighborhood of $S_k$, there are unique $y(x) \in S_k$ and $z(x) \geq 0$ such $x = y(x) + z(x) n(y(x))$.
Let $\varphi(x) = (y(x), z(x))$; for all sufficiently small $\delta$, this is a bijection between a neighborhood $S_k^\delta$ of $S_k$ 
and the set $S_k \times [0, \delta)$, i.e., $(y,z)$ can be viewed as a new set of coordinates on  $S_k^\delta$. 

Define 
\[
h_1 =  L z,~~h_2 = \frac{1}{2} Lz^2.
\]
From the assumptions on the vector fields $v_0,...,v_n$, it follows (see Lemma~\ref{gene}) that there are $\alpha_k, \beta_k \in C^1(S_k)$such that $\alpha_k(y) \geq 0$ for $y \in S_k$ and 
\begin{equation} \label{abeta} 
h_1(x) = \beta_k(y(x)) z(x) + O(z^2(x)),~~h_2(x) = \alpha_k(y(x)) z^2(x) + O(z^3(x)),~~{\rm as}~~z(x) \rightarrow 0. 
\end{equation}
The functions $\beta_k$ and $\alpha_k$ are the leading terms for the drift and diffusion coefficients, respectively, in the direction 
co-normal (with respect to $L^\varepsilon$)  
to the boundary if the operator $L$ is written in $(y,z)$ coordinates.
We will assume that $\alpha_k > 0$ for each $y \in S_k$ (we could weaken this assumption and instead assume that there is $y \in S_k$ such that $\alpha_k(y) > 0$).  
Since $X^x_t$ is a non-degenerate diffusion on $S_k$, it has a unique invariant probability measure on $S_k$. This measure will be denoted by $\pi_k$. 
Define 
\[
\bar{\alpha}_k = \int_{S_k} \alpha_k(y) d\pi_k(y),~~\bar{\beta}_k = \int_{S_k} \beta_k(y) d\pi_k(y) .
\]
 We will see that  if $\bar{\alpha}_k > \bar{\beta}_k $,
then $\mathrm{P} (\lim_{t \rightarrow \infty} {\rm dist}(X^x_t, S_k) = 0) > 0$ for each $x \in \overline{D}$ in a sufficiently small neighborhood of $S_k$. If $\bar{\alpha}_k < \bar{\beta}_k$, then this
probability is zero unless $x \in S_k$. We will assume that $\bar{\alpha}_k \neq \bar{\beta}_k$ for each $k$  and will
refer to $S_k$ as attracting if $\bar{\alpha}_k > \bar{\beta}_k $ and repelling if $\bar{\alpha}_k < \bar{\beta}_k$.
\\

In this section and the next one, we will examine the behavior of $X^{x, \varepsilon}_t$ near $S_k$. Since we assume that $k$ is fixed, we will temporarily drop this subscript from the notation. From the conditions placed on the vector fields $v_0,v_1,...,v_d$, it follows that the operator $L$ can be applied to functions defined on $S = S_k$. 
In order to stress that we are considering the restriction of the operator to $S$ (where variables $y$ are used), we denote the resulting operator by $L_y$.  Note that the operator
$L_y$ acting in the $y$ variables can be applied to functions of $(y,z)$ by treating $z$ as a parameter. 
\begin{lemma} \label{gene} The generator of the process $X^x_t$ in $(y,z)$ coordinates can be written as:
\begin{equation} \label{opr}
L u  = L_y u  +  z^2 \alpha(y) \frac{\partial^2 u}{\partial z^2} +   z  \beta(y) \frac{\partial u}{\partial z}  + z {\mathcal{D}}_y \frac{\partial u}{\partial z} + R u 
\end{equation}
with
\[
R u = z {\mathcal{K}}_y u   + z^2 \mathcal{N}_y \frac{\partial u}{\partial z}  +  z^3 \sigma(y,z) \frac{\partial^2 u}{\partial z^2},
\]
where ${\mathcal{D}}_y$ is a first-order differential operator in $y$ (without a zero-order term), whose coefficients depend only on the $y$ variables, ${\mathcal{K}}_y$ is a  differential operator  on $S \times [0, \delta)$ with first- and second-order derivatives in $y$,  ${\mathcal{N}}_y$ is a differential operator on 
$S \times [0, \delta)$ with first-order derivatives in $y$ and a potential term. All the 
operators have $C^1$ coefficients,  
while $\alpha, \beta \in C^1(S), \sigma \in C^1(S^\delta)$. 
The generator of the process $X^{x,\varepsilon}_t$ in $(y,z)$ coordinates is the operator
\[
L^\varepsilon = L + \varepsilon^2   \tilde{L},
\]
where $\tilde{L}$ is a second-order uniformly elliptic differential operator with $C^1$ coefficients in $(y,z)$ variables. 
\end{lemma}
\proof
Consider the simplest case: $S$ is one-dimensional (i.e., $d = 2$) and defined, in $(y,z)$ coordinates, in a 
neighborhood of a point $(y_0, z_0) \in S$, as the line $\{(y,z): z = 0\}$. (The case when $d > 2$ requires only slightly more complicated notations, while
the case when $S$ is a curve (surface) can be reduced to the case when it is a linear subspace by a change of variables.) Each term of the operator $L$ defined in (\ref{opl}) can
be considered separately, so we can assume that
\[
Lu = 
\frac{1}{2} \frac{\partial}{\partial v}(\frac{\partial u}{\partial v}),
\] 
where $v = (v^1, v^2)$ is a vector field tangent to $S$ (the term with the first order derivative can be considered similarly). Thus
\[
Lu = 
\frac{1}{2}(v^1)^2 \frac{\partial^2 u}{\partial y^2} +  v^1 v^2  \frac{\partial^2 u}{\partial y \partial z} + 
\]
\[
\frac{1}{2} (v^2)^2 \frac{\partial^2 u}{\partial z^2} + 
\frac{1}{2} (v^1 \frac{\partial v^1}{\partial y} + v^2 \frac{\partial v^1}{\partial z}) \frac{\partial u}{\partial y} + \frac{1}{2} (v^1 \frac{\partial v^2}{\partial y} +
v^2 \frac{\partial v^2}{\partial z}) \frac{\partial u}{\partial z}.
\]
Using the smoothness of $v^1, v^2$ and the fact that $v^2(y,0) = 0$, we can write
\[
v^1(y,z) = v^1(y,0) + z g_1(y,z),~~v^2(y,z) = \frac{\partial v^2}{\partial z} (y,0) z + z^2 g_2(y,z),~~
\]
\[
\frac{\partial v^1}{\partial y}(y,z) = \frac{\partial v^1}{\partial y}(y,0) + z g_3(y,z),~~ 
\frac{\partial v^1}{\partial z}(y,z) = \frac{\partial v^1}{\partial z}(y,0) + z g_4(y,z),
\]
\[
\frac{\partial v^2}{\partial y}(y,z) = \frac{\partial^2 v^2}{\partial y \partial z} (y,0) z + z^2 g_5(y,z),~~\frac{\partial v^2}{\partial z}(y,z) = 
\frac{\partial v^2}{\partial z}(y,0) + z g_6(y,z),
\]
where $g_1,...,g_6$ are smooth functions. Expressing the coefficients of $L$ using these expansions, we obtain the desired form of the operator. 
The statement about the form of $L^\varepsilon$ follows immediately.
\qed
\\

\noindent
{REMARK 1.} From Lemma~\ref{gene} it follows that the functions $\alpha$ and $\beta$ found in the coefficients of $L$ agree with those ($\alpha_k$ and $\beta_k$)
defined in (\ref{abeta}). Moreover, $\alpha$ and $\beta$ do not depend on the choice of  $n$, i.e., we could
replace the field of co-normals by any sufficiently smooth field of vectors forming acute angles with inward-pointing normals. 
\\
\\
{REMARK 2.} Lemma~\ref{gene} shows that the operator $L$  is approximated by the operator
\begin{equation} \label{kok}
K u  = L_y u  +  z^2 \alpha(y) \frac{\partial^2 u}{\partial z^2} +   z  \beta(y) \frac{\partial u}{\partial z}  + 
z {\mathcal{D}}_y \frac{\partial u}{\partial z}.
\end{equation}
The latter operator has a useful homogeneity property: $(Ku)(y,az) = K (u(y,az))$ for each $y$ and $a >0$. 
\\

The next lemma will be instrumental in analyzing the time it takes the process $X^{x, \varepsilon}_t$ to
leave a small neighborhood of $S$ if it starts at $x \in S$ (in the attracting case) or to reach $S$ if it starts at a point close to $S$ (in the repelling case).  
\begin{lemma} \label{spec}
If $\bar{\alpha} > \bar{\beta}$ ($\bar{\alpha} < \bar{\beta}$), then there exist $\gamma >0$ ($\gamma <0$) and a positive-valued function 
$\varphi  \in C^1(S)$ satisfying $\int_S \varphi d \pi = 1$ such that
\begin{equation} \label{mg}
L_y \varphi + \alpha \gamma (\gamma -1) \varphi + \beta \gamma \varphi + \gamma \mathcal{D}_y \varphi = 0.
\end{equation}
Such $\gamma$ are $\varphi$ are determined uniquely. 
\end{lemma}
\proof Let $\lambda_\gamma$ be the top eigenvalue for $M(\gamma)$, where $M(\gamma)$ is the operator in the left hand side of (\ref{mg}).
Since $L_y + \gamma \mathcal{D}_y$ is an elliptic operator on a compact manifold, this eigenvalue is simple, as follows from the 
Perron-Frobenius Theorem by considering the corresponding parabolic semigroup whose time one kernel is positive. 
Therefore  $\lambda_\gamma$ depends smoothly on the parameter $\gamma$  (\cite{Kato}). Moreover, the 
 corresponding eigenfunction $\varphi_\gamma$ can be chosen so that it depends smoothly on $\gamma$ and $\varphi_0 \equiv 1$. 
Differentiating the equality
$ M(\gamma) \varphi_\gamma = \lambda_\gamma \varphi_\gamma $ in $\gamma$, we obtain
\[
M(\gamma) \varphi'_\gamma + (\alpha (2 \gamma - 1) + \beta +  \mathcal{D}_y) \varphi_\gamma = \lambda'_\gamma \varphi_\gamma +  \lambda_\gamma \varphi'_\gamma.
\]
Put $\gamma = 0$ and integrate both sides with respect to the measure $\pi$. 
Note that $\int_S \varphi_0 d \pi = 1$, $\mathcal{D}_y \varphi_0 = 0$, and $\lambda_0 = 0$.
Since $M(0) = L_y$ and $\pi$ is invariant for the process generated by $L_y$, we have
$\int_S M(0) \varphi'_0 d \pi  = 0$. Thus 
\[
 \lambda'_0  = \int_S (-\alpha + \beta) d \pi = \bar{\beta} - \bar{\alpha}.
\]
Let us assume that $\bar{\beta} <  \bar{\alpha}$ (the case when $\bar{\beta} >  \bar{\alpha}$ can be handled similarly). 
Then $ \lambda'_0 < 0$. Observe that $\lambda_0 = 0$. 
Let us demonstrate that
  $\lim_{\gamma \rightarrow \infty} \lambda_\gamma = +\infty$. Indeed, let $\pi_\gamma$ be the invariant probability measure for the
process governed by the operator $L_y + \gamma \mathcal{D}_y$ on $S$. Integrating the equality 
$ M(\gamma) \varphi_\gamma = \lambda_\gamma \varphi_\gamma $ with respect to $\pi_\gamma$, and dividing both sided by $\gamma$, we obtain
\[
(\gamma - 1) \int_S \alpha \varphi_\gamma d \pi_\gamma + \int_S \beta \varphi_\gamma d \pi_\gamma = \gamma^{-1} \lambda_\gamma \int_S \varphi_\gamma 
d \pi_\gamma.
\]
Since $\alpha > 0$ and $\beta$ is bounded, this implies that $\lim_{\gamma \rightarrow \infty} \lambda_\gamma = +\infty$. Therefore, there exists $\gamma > 0$ such that $\lambda_\gamma = 0$. 

Let us 
show that such $\gamma$ is unique.  Assume the contrary, i.e., that $(\gamma_1, \varphi_1)$ and $(\gamma_2, \varphi_2)$ satisfy 
(\ref{mg}) and $0 < \gamma_1 < \gamma_2$.  (The case when $\gamma_1$ and $\gamma_2$ are negative can
be considered similarly.)  Let $\tilde{X}^x_t 
= (\tilde{Y}^x_t, \tilde{Z}^x_t)$ be the
family of diffusion processes on $S \times (0, \infty)$ with the generator $K$ defined in (\ref{kok}).
(This operator can be obtained from the generator of  ${X}^x_t$ by discarding the last term in~(\ref{opr}).) 

Consider the processes $\xi^1_t =  \varphi_1(\tilde{Y}^x_t)(\tilde{Z}^x_t)^{\gamma_1}$ and 
$\xi^2_t = \varphi_2(\tilde{Y}^x_t)(\tilde{Z}^x_t)^{\gamma_2}$. Here, we fix an initial point $x$ such that $\xi^1_0 = 1$. Let 
\[
\tau_n = \inf\{t: \xi^1_t = \frac{1}{n}~{\rm or}~ \xi^1_t = n\}.
\]
Since $K (\varphi_1(y) z^{\gamma_1}) = 
K (\varphi_2(y) z^{\gamma_2}) = 0$,  by the Ito formula,
$\xi^1_t$ and $\xi^1_t$ are local martingales, while the stopped processes $\xi^1_{\tau_n \wedge t}$ and $\xi^2_{\tau_n \wedge t}$ are martingales since they are bounded. Note that $\tau_n < \infty$ almost surely since $\gamma_1 
\neq 0$. By the 
Optional Stopping Theorem applied to the process $\xi^1_{\tau_n \wedge t}$, $\mathrm{E} \xi^1_{\tau_n} = \xi^1_0 = 1$, and
therefore $\mathrm{P}(\xi^1_{\tau_n} = 1/n) = n/(n+1)$, while 
$\mathrm{P}(\xi^1_{\tau_n} = n) = 1/(n+1)$.  We estimate
\[
\mathrm{E} \xi^2_{\tau_n} \geq \inf ({\varphi_2}/{\varphi_1^{\frac{\gamma_2}{\gamma_1}}})  \mathrm{E} 
( \xi^1_{\tau_n})^{\frac{\gamma_2}{\gamma_1}} \geq   \inf ({\varphi_2}/{\varphi_1^{\frac{\gamma_2}{\gamma_1}}}) \mathrm{P}(\xi^1_{\tau_n} = 
n) n^{\frac{\gamma_2}{\gamma_1}} \rightarrow \infty~~{\rm as}~n \rightarrow \infty,
\]
since $\mathrm{P}(\xi^1_{\tau_n} = 
n) n^{{\gamma_2}/{\gamma_1}} = (n+1)^{-1} n^{{\gamma_2}/{\gamma_1}} \rightarrow \infty$.
However, by the 
Optional Stopping Theorem applied to the process $\xi^2_{\tau_n \wedge t}$, $\mathrm{E} \xi^2_{\tau_n} = \xi^2_0$ does not depend on $n$.
Thus we get a contradiction, which proves uniqueness. \qed
\\

\noindent
{REMARK.} If $\beta/\alpha = {\rm const}$, then $\gamma = 1 - \beta/\alpha$ and $\varphi \equiv 1$ solve (\ref{mg}).

\section{Behavior of the process $X^{x,\varepsilon}_t$ near the boundary} \label{neigh}
For a closed set $A$, let $\tau^{x,\varepsilon}(A) = \inf\{t \geq 0: X^{x,\varepsilon}_t \in A\}$.
 In this section, as before, we drop the 
subscript $k$ from the notation since we are going to focus on a small neighborhood of a single surface $S = S_k$. 
 For sufficiently small $\varkappa \geq 0$, we can define 
\[
\Gamma_\varkappa = \{(y,z): (\varphi(y))^{\frac{1}{\gamma}} z = \varkappa\},
\]
where $\varphi$ is defined in Lemma~\ref{spec}. For $\varkappa_1< \varkappa_2$ that are sufficiently small, we denote the region
between $\Gamma_{\varkappa_1}$ and $\Gamma_{\varkappa_2}$ by
\[
V_{\varkappa_1, \varkappa_2} = \{(y,z): \varkappa_1 \leq (\varphi(y))^{\frac{1}{\gamma}} z  \leq  \varkappa_2\}.
\]
In this section, we will study the asymptotics, as $\varepsilon \downarrow 0$, of the following quantities:

(a) The measure on $S$ that is induced by 
$X^{x,\varepsilon}_{ \tau^{x,\varepsilon} (S)}$, where $x \in \Gamma_\varkappa$. 

(b) The measure on $\Gamma_\varkappa$ that is induced by 
$X^{x,\varepsilon}_{ \tau^{x,\varepsilon} (\Gamma_\varkappa)}$, where $x \in S$. Here, $\varkappa > 0$ is assumed to be
sufficiently small. 

(c) The stopping time $\tau^{x,\varepsilon} (\Gamma_\varkappa)$ 
(particularly $\mathrm{E} \tau^{x,\varepsilon} (\Gamma_\varkappa)$) for $x \in S$. Here, $\varkappa > 0$ is assumed to be
sufficiently small. 

(d) The stopping time $\tau^{x,\varepsilon} (S \bigcup \Gamma_{\varkappa})$ 
(particularly $\mathrm{E} \tau^{x,\varepsilon} (S \bigcup \Gamma_{\varkappa})$) for 
$x \in \Gamma_\zeta$. Here, $0 < \zeta <  \varkappa$ and   $\varkappa$ is assumed to be
sufficiently small.


\subsection{The measure induced by the process stopped at the boundary} \label{sseca}

The asymptotics, as $\varepsilon \downarrow 0$, of the measure on $S$   induced by 
$X^{x,\varepsilon}_{ \tau^{x,\varepsilon} (S)}$,   $x \in \Gamma_\varkappa$, was studied in
\cite{FK21}. Here, we recall the main arguments since some of
these will be needed in other parts of the paper, and also to make the presentation self-contained.

Formula (\ref{opr}) shows that the generator of $X^{x}_t$ can be approximated near the boundary by an operator that is homogeneous in the variable $z$. The  remainder term with the operator $R$ can be made small by considering a sufficiently small neighborhood of $S$. The same approximation is useful for the generator of $X^{x, \varepsilon}_t$, except in a yet smaller neighborhood of $S$, where the 
perturbation $\varepsilon^2   \tilde{L}$ is comparable to or larger than the operator $L$, since the coefficients of the latter degenerate near $S$. In order to understand the behavior of $X^{x, \varepsilon}_t$ in an $\varepsilon$-dependent neighborhood of $S$, 
we introduce the change of variables
\[
\Psi_\varepsilon(y,z) = (y,\frac{z}{\varepsilon}),~~~\Psi_\varepsilon: S \times[0,\delta) \rightarrow S \times [0,\frac{\delta}{\varepsilon})
\subset S \times [0,\infty)
\]
and the operator $M^\varepsilon u  =  {L}^\varepsilon (u (\Psi_\varepsilon)) (\Psi_\varepsilon^{-1})$.
This operator is the generator of the  process 
\begin{equation} \label{chp}
\mathcal{X}^{{\bf x},\varepsilon}_t := \Psi_\varepsilon (X^{\Psi_\varepsilon^{-1} ({\bf x}), \varepsilon}_t)
\end{equation}
 on $S \times [0,{\delta}/{\varepsilon})$.  To stress the difference between the two sets of coordinates, we use the notation
${\bf x} = (y,{\bf z})$ for the new variables instead of $x = (y,z)$. 
By Lemma~\ref{gene},  
\begin{equation} \label{genpp0}
M^\varepsilon u  = L_y u  +  ({\bf z}^2 \alpha(y) + \rho(y))\frac{\partial^2 u}{\partial {\bf z}^2} +   
{\bf z} \beta(y) \frac{\partial u}{\partial {\bf z}}  + 
{\bf z} {\mathcal{D}}_y \frac{\partial u}{\partial {\bf z}} + 
\widehat{R}^\varepsilon u =: M u +  \widehat{R}^\varepsilon u, 
\end{equation}
where $\rho$ is the coefficient in front of the second derivative in the variable $z$ at $z=0$ in the operator $\tilde{L}$, and $\widehat{R}^\varepsilon$ is a second order operator with continuously differentiable coefficients 
that tend to zero uniformly on $S \times [0,r]$ as $\varepsilon \downarrow 0$ for each $r > 0$. Thus, in a small neighborhood of $S$, $M^\varepsilon$ is a small perturbation of the operator $M$, which does not depend on $\varepsilon$. 

We can view $M$ as an operator on $S \times [0,\infty)$.
Let us examine the behavior of the process with the generator $M$ (which differs from a homogeneous operator by the
presence of the extra term $\rho(y){\partial^2 u}/{\partial {\bf z}^2}$). The process with the generator $M$, earlier denoted by
$\mathcal{{X}}^{\bf x}_t$, will also be written as 
$(\mathcal{Y}^{{\bf x}}_t, \mathcal{Z}^{{\bf x}}_t)$, where ${\bf x}=(y,{\bf z})$ is the initial point; its state space is $S \times [0,\infty)$;  the
process is stopped upon reaching $S \times \{0\}$. 
Let  ${\eta}^{\bf x} = \inf\{t \geq 0: \mathcal{Z}^{\bf x}_t  = 0\}$ for ${\bf x} \in S \times [0, \infty)$. 

\begin{lemma} \label{sttm}
If  $\bar{\alpha} > \bar{\beta} $ (the boundary is attracting), then $\mathrm{P}({\eta}^{\bf x} < \infty) = 1$ for each ${\bf x} = (y,{\bf z}) \in S \times [0, \infty)$. 
If $\bar{\alpha} < \bar{\beta} $ (the boundary is repelling), 
then $\lim_{{\bf z} \rightarrow \infty} \mathrm{P}({\eta}^{(y,{\bf z})} < \infty) = 0$ uniformly in $y \in S$. 
\end{lemma}
\proof Let $\Phi(y,{\bf z}) = (y,\ln({\bf z}))$ be the mapping from $S \times (0,\infty)$ to $S \times\mathbb{R}$.
Consider the process $\Phi(\mathcal{{X}}^{\bf x}_t) = (\mathcal{Y}^{\bf x}_t, \ln(\mathcal{Z}^{\bf x}_t)) $ 
on $S \times \mathbb{R}$. This process may go to $-\infty$ along the ${\bf z}$-axis in finite time, but this will not cause any problems. The generator of this process is 
\begin{equation} \label{genaa}
\mathcal{A} u   = L_y u  +  (\alpha(y) + \rho(y) e^{-2 {\bf z}})(\frac{\partial^2 u}{\partial  {\bf z}^2} - 
\frac{\partial u}{\partial {\bf z}} ) +
\beta(y) \frac{\partial u}{\partial  {\bf z}}  + 
{\mathcal{D}}_y \frac{\partial u}{\partial  {\bf z}}. 
\end{equation}
Let $\psi: S \rightarrow \mathbb{R}$ solve
\[
L_y \psi = \alpha - \beta - (\bar{\alpha} - \bar{\beta}),~~~\int_S \psi d \pi  = 0.
\]
Note that the solution of $L_y \psi = \alpha - \beta - (\bar{\alpha} - \bar{\beta})$ exists since $\int_S  (\alpha - \beta - (\bar{\alpha} - \bar{\beta})) d \pi = 0$, and is defined uniquely up to the solutions of the homogeneous equation (i.e., up to a constant). Therefore, the 
condition that the solution is orthogonal to $\pi$ determines $\psi$ uniquely. 
 Let $g(y,{\bf z}) =    \psi(y) +  {\bf z} $.
  From the Ito formula applied to
$g(\mathcal{Y}^{\bf x}_t,  \ln(\mathcal{Z}^{\bf x}_t) )$, it follows that
\[
h^{{\bf x}}_t : =  \psi(\mathcal{Y}^{{\bf x}}_t)  + 
\ln(\mathcal{Z}^{\bf x}_t) - \int_0^t \mathcal{A} g(\mathcal{Y}^{{\bf x}}_s, \ln(\mathcal{Z}^{{\bf x}}_s)) ds = 
\]
\[
 \psi(\mathcal{Y}^{{\bf x}}_t) + 
\ln(\mathcal{Z}^{\bf x}_t) +  \int_0^t \rho(\mathcal{Y}^{\bf x}_s) (\mathcal{Z}^{{\bf x}}_s)^{-2} ds + (\bar{\alpha} - \bar{\beta}) t
\]
is a local martingale. For $n \in \mathbb{Z}$, define the following subsets of $ S \times (0,\infty)$:
\begin{equation} \label{subsets}
  B_n     = \{{\bf x}: g( \Phi({\bf x}) ) = n\},~~~~ B_n^- = \{{\bf x}: g( \Phi({\bf x}) ) \leq n\},~~~~ B_n^+ =
	\{{\bf x}: g( \Phi({\bf x}) ) \geq n\}.
\end{equation}
Since $({\mathcal{Y}}^{\bf x}_t, {\mathcal{Z}}^{\bf x}_t)$ is a non-degenerate diffusion, there is a constant $c > 0$ such that
\begin{equation} \label{recc}
\mathrm{P} ({\eta}^{\bf x} < \infty) > c,~~{\bf x} \in  B_0^-.
\end{equation}
 For $\bar{\alpha} > \bar{\beta} $, 
the process  $\psi(\mathcal{Y}^{{\bf x}}_t) + 
\ln(\mathcal{Z}^{\bf x}_t) $ is a  local supermartingale. Since it is unbounded with probability one, the
process $\mathcal{{X}}^{\bf x}_t = (\mathcal{Y}^{\bf x}_t, \mathcal{Z}^{{\bf x}}_t)$ reaches $B_0^-$ with probability one for each initial point ${\bf x}$. 
From (\ref{recc}) and the strong Markov property, it follows that $\mathrm{P} ({\eta}^{\bf x} < \infty) = 1$ for each ${\bf x}$.  

Now assume that $\bar{\alpha} < \bar{\beta} $. For ${\bf x} \in B_n$, let 
\[
\sigma^{{\bf x}} = \inf\{t \geq 0: \mathcal{{X}}^{\bf x}_t \in B_{n-1} \bigcup B_{n+1}\}.
\]
It is sufficient to show that there exists $c > 0$ such that, for all sufficiently large~$n$,
\begin{equation} \label{slnn}
\mathrm{P} ( \mathcal{{X}}^{\bf x}_{\sigma^{{\bf x}} }  \in B_{n+1}) \geq \frac{1}{2} + c,~~~
{\bf x} \in B_n.
\end{equation}
Since $(\mathcal{Y}^{\bf x}_t, \ln(\mathcal{Z}^{\bf x}_t))$ is a diffusion with coefficients that are bounded on $B^+_0$, 
there is $c' > 0$ such
that, 
\[
\mathrm{E} \sigma^{{\bf x}} \geq c',~~~
{\bf x} \in B_n,
\]
provided that $n \geq 1$. Since $\mathrm{E}(h^{{\bf x}}_{\sigma^{{\bf x}}} - h^{{\bf x}}_0) = 0$,
\[
2 \mathrm{P} ( \mathcal{{X}}^{\bf x}_{\sigma^{\bf x}}  \in B_{n+1}) - 1 
- \mathrm{E} \int_0^{ \sigma^{{\bf x}}} \left( (\bar{\beta} - \bar{\alpha} ) - \rho(\mathcal{Y}^{\bf x}_s) (\mathcal{Z}^{{\bf x}}_s)^{-2}
\right) ds  = 0.
\]
For all sufficiently large $n$, the integrand in the last integral is estimated from below by $(\bar{\beta} - \bar{\alpha})/2$.
Therefore,
\[
\mathrm{P} ( \mathcal{{X}}^{\bf x}_{\sigma^{\bf x}} \in B_{n+1}) \geq \frac{1}{2} + 
\frac{c' (\bar{\beta} - \bar{\alpha})}{4},~~~
{\bf x} \in B_n,
\]
as required.
\qed
\\

 The solutions to the equations in the following lemma
are sought in the spaces $\mathcal{C} : = C^2(S \times (0,\infty)) \bigcap C_b(S \times [0,\infty))$ or 
$\mathcal{C}^2 := C^2(S \times [0,\infty))$. 
\begin{lemma} \label{mleee}
Let $f \in C(S)$.

(a)  If the boundary $S$ is attracting, then there is a unique solution $u \in \mathcal{C}$ to the equation
\begin{equation} \label{fieq}
M u (y,{\bf z}) = 0,~~y \in S, {\bf z} > 0;~~~~u(y,0) = f(y).
\end{equation}

(b) If the boundary $S$ is repelling, then there is a unique solution $h \in \mathcal{C}^2$ to the equation
\begin{equation} \label{fieq2}
M h (y,{\bf z}) = 0,~~y \in S, {\bf z} > 0;~~~~h(y,0) \equiv 1;~~~\lim_{{\bf z} \rightarrow \infty} \sup_{y \in S} |h(y,{\bf z})|  = 0.
\end{equation}
There is a unique solution  $u \in \mathcal{C}$ to the equation
\begin{equation} \label{fieq3}
M ( h u) (y,{\bf z}) = 0,~~y \in S, {\bf z} > 0;~~~~u(y,0) = f(y).
\end{equation}

(c) In both cases, there exists a constant $\overline{u}$ such that  $\overline{u}= \lim_{{\bf z} \rightarrow \infty} u (y,{\bf z})$, uniformly in $y \in S$.
\end{lemma}

\proof If $S$ is attracting, we define $u({\bf x}) = \mathrm{E} f(\mathcal{X}^{\bf x}_{{\eta}^{\bf x}})$,
where the right-hand side is well-defined since $\mathrm{P}({\eta}^{\bf x} < \infty) = 1$ (Lemma~\ref{sttm}). This is a standard
probabilistic representation of the solution to equation (\ref{fieq}), and the solution is unique in~$\mathcal{C}$ (see, e.g., \cite{F85}).

If $S$ is repelling, we define $h({\bf x}) = \mathrm{P} ({{\eta}^{\bf x}} < \infty)$. By Lemma~\ref{sttm},
$\lim_{{\bf z} \rightarrow \infty} \mathrm{P}({\eta}^{(y,{\bf z})} < \infty) = 0$, and thus $h \in  \mathcal{C}^2$ is 
the unique solution to (\ref{fieq2}).  Moreover, the process ${\widehat{\mathcal{X}}}^{\bf x}_t$, defined by conditioning 
${\mathcal{X}}^{\bf x}_t$ on the event 
$\{{{\eta}^{\bf x}} < \infty\}$, is governed by the operator $\widehat{M} u = h^{-1} M ( h u)$ (the Doob transform of $M$, 
see, e.g., \cite{Pin}). The operator is non-degenerate, 
and, by construction, the process reaches $S$ with probability one for each initial point ${\bf x}$ (see, e.g., \cite{Pin}). Therefore, there is a unique solution
$u \in \mathcal{C}$ to (\ref{fieq3}), which is given by $u({\bf x}) = 
\mathrm{E} f(\widehat{\mathcal{X}}^{\bf x}_{\widehat{\eta}^{\bf x}})$,  where
$ \widehat{\eta}^{\bf x} = \inf\{t \geq 0: \widehat{\mathcal{X}}^{\bf x}_t \in S \times \{0\}\} $.

It remains to prove part (c) of the lemma. We will need the following fact. Suppose that the generator 
of a diffusion process $H^x_t$ is a uniformly elliptic operator in a bounded domain $G$ with a smooth boundary. Let $K \subset G$ 
be compact. Let $\tau^x = \inf\{t \geq 0: H^x_t \in \partial G$. Let   $\mu^x$ be the measure on $\partial G$ induced by $H^x_{\tau^x}$,
and let $p^x$ be its density with respect to the Lebesgue measure. Then there is a constant $c > 0$ such that
\begin{equation} \label{esbl}
p^x(\tilde{x}) \geq c,~~x \in K, \tilde{x} \in \partial G.
\end{equation}
The bound $c$ can be chosen to be the same for all the operators that have the same ellipticity constant and bound on the $C$-norm of the coefficients. The estimate (\ref{esbl}) is likely known in the PDE literature, but since we could not find a proof, we 
sketched it in \cite{FK21} using an estimate on the Poisson kernel for the generator of $H^x_t$ in a small inward-directed cone with 
the vertex at $\tilde{x}$ (see \cite{Va}).

First, consider the case when $S$ is attracting. 
Recall the definition of the sets $B_n$ and $B_n^+$ from (\ref{subsets}). For ${\bf x} \in B_n^+$, let
\[
\sigma_n^{{\bf x}} = \inf\{t \geq 0: \mathcal{{X}}^{\bf x}_t \in B_{n}\}.
\]
Then, for ${\bf x} \in B_n^+$,
\begin{equation} \label{fkf}
u({\bf x}) = \mathrm{E} u( \mathcal{{X}}^{\bf x}_{\sigma_n^{{\bf x}} } ).
\end{equation}
Let $V_n = \sup_{{\bf x}_1,{\bf x}_2 \in B_n}|u({\bf x}_1) - u({\bf x}_2)|$. Then $V_0 \leq  \sup_{x_1,x_2 \in S}|f(x_1) - f(x_2)|$ and,
by (\ref{fkf}), 
\[
\sup_{{\bf x}_1,{\bf x}_2 \in B^+_n}|u({\bf x}_1) - u({\bf x}_2)| \leq V_n.
\]
Thus it is sufficient to show that $V_n \rightarrow 0$ as $n \rightarrow \infty$. Since the operator $\mathcal{A}$ in (\ref{genaa}) is uniformly elliptic and
its coefficients are bounded (uniformly in $n \geq 0$) in the domain bounded by $\Phi(B_{n})$ and $\Phi(B_{n+2})$, (\ref{esbl}) 
is applicable to the process $\Phi(\mathcal{{X}}^{\bf x}_t)$  with  $K = \Phi(B_{n+1})$. Consequently, for the density 
$p^{\bf x}_n$ of the measure $\mu^{\bf x}_n$ induced by 
$ \Phi(\mathcal{{X}}^{\bf x}_{\sigma_n^{{\bf x}} })$ on $\Phi(B_n)$, we have 
$p^{\bf x}_n(\tilde{{\bf  x}}) \geq c$, 
${\bf x} \in B_{n+1}$, $\tilde{{\bf x}} \in  \Phi(B_n)$. (Here we used the fact that $p^{\bf x}_n \geq \overline{p}^{\bf x}_n$, where
$\overline{p}^{\bf x}_n$ corresponds to stopping the process $\mathcal{{X}}^{\bf x}_t$ 
on $B_n \bigcup B_{n+2}$ rather than on $B_n$.)
Therefore, by (\ref{fkf}), for ${\bf x}_1,{\bf x}_2 \in B_{n+1}$,
\[
|u({\bf x}_1)  -u({\bf x}_2)| \leq (1 - c \lambda(\Phi(B_n))) V_n,
\]
where $\lambda(\Phi(B_n)) = \lambda(\Phi(B_0))$ is the Lebesgue measure of $\Phi(B_n)$. Thus, $V_n \leq V_0 (1 - c 
\lambda(\Phi(B_0)))^n \rightarrow 0$ as $n \rightarrow \infty$, which implies that there is a limit
\[
\overline{u} = \lim_{{\bf z} \rightarrow \infty} u (y,{\bf z}) 
\]
uniformly in $y \in S$. 

The same argument applies in the case  when $S$ is repelling. We only need to observe that the operator governing the 
conditioned process $ \Phi(\widehat{\mathcal{X}}^{\bf x}_t)$ 
in the domain  between $\Phi(B_n)$ and $\Phi(B_{n+2})$ has coefficients bounded from above and the ellipticity constant bounded from below uniformly in $n \geq 0$,  as follows from the standard 
elliptic estimates on the function $h$. \qed
\\

In the proof of Lemma~\ref{mleee}, we saw that $u({\bf x}) = \mathrm{E} f(\mathcal{X}^{\bf x}_{{\eta}^{\bf x}})$ if $S$ is
attracting, and $u({\bf x}) = \mathrm{E} f(\widehat{\mathcal{X}}^{\bf x}_{\widehat{\eta}^{\bf x}})$ if $S$ is repelling. 
Let $\nu^{\bf x}$ be the measure on $S$ induced by $\mathcal{X}^{\bf x}_{{\eta}^{\bf x}}$ in the former case, or by 
$\widehat{\mathcal{X}}^{\bf x}_{\widehat{\eta}^{\bf x}}$ in the latter case. Thus $u({\bf x}) = \int_S f d \nu^{\bf x}$. Since the 
mapping $f \rightarrow \overline{u}$ is linear and continuous from $C(S)$ to $\mathbb{R}$, there is a measure $\nu$ on $S$ such that
$\overline{u} = \int_S f d \nu$. From part (c) of Lemma~\ref{mleee}, it follows that $\nu^{(y, {\bf z})} \rightarrow \nu$ weakly, as
${\bf z} \rightarrow \infty$.



\begin{lemma} \label{mlee1}
Let $f \in C(S)$.  For each $\varkappa > 0$ (such that $\Gamma_\varkappa$ is defined),
\[
\lim_{\varepsilon \downarrow 0} \int_S f d \nu^{x,\varepsilon}  = \int_S f d \nu
\]
uniformly in $x \in \Gamma_\varkappa$, where $\nu^{x,\varepsilon}$ is the measure on $S$ induced by 
$X^{x,\varepsilon}_{ \tau^{x,\varepsilon} (S)}$. 
\end{lemma}

\proof Recall that $\mathcal{X}^{\bf x}_t$ is the process on $S \times [0,\infty)$ with the generator 
$M$ and $\mathcal{X}^{{\bf x},\varepsilon}_t = ({\mathcal{Y}}^{{\bf x},\varepsilon}_t,{\mathcal{Z}}^{{\bf x},\varepsilon}_t)$ is 
the process with the generator $M^\varepsilon$. For each $r > 0$, the latter operator can be defined on $S \times [0,r]$, 
provided that $\varepsilon > 0$  sufficiently small. 

 Let  $\delta > 0$.   
Recall that $u({\bf x}) = \mathrm{E} f(\mathcal{X}^{\bf  x}_{{\eta}^{\bf x}})$
is the function from Lemma~\ref{mleee}, $\overline{u} = \int_S f d \nu$, and recall the surfaces $B_n$ defined in (\ref{subsets}).
Let $u^\varepsilon(x) = \mathrm{E} f(X^{x,\varepsilon}_{ \tau^{x,\varepsilon} (S)})$.  We will also write $u^\varepsilon({\bf x}) $ 
when the same function is considered in $(y,{\bf z})$ coordinates. 
By the strong Markov property of the process, it is sufficient to show that 
\[
|u^\varepsilon({\bf x}) -  \overline{u}|   \leq \delta, ~~~{\bf x} \in B_n.
\]
First, consider the case when  $S$ is attracting. Let $n \in \mathbb{N}$ be such that $|u({\bf x}) - \overline{u}| \leq \delta/2$ for
${\bf x} \in B_n$.  
Fix $\eta > 0$, to be specified later. Take $T, r > 0$ such that
\begin{equation} \label{bnn}
\mathrm{P}({\eta}^{\bf x} > T) \leq \eta,~~~ \mathrm{P}(\sup_{0 \leq t \leq {\eta}^{\bf x}} {\mathcal{Z}}^{\bf x}_t > r-1)  \leq \eta,~~~
{\bf x}\in B_n.
\end{equation}
Such $T$ and $r$ exist  since $\mathrm{P}({\eta}^{\bf x} < \infty) = 1$ and the probabilities in the left-hand side of
both inequalities depend continuously on ${\bf x}$. Define
\[
{\eta}^{{\bf  x},\varepsilon}_r =  \begin{cases} \inf\{t \geq 0: {\mathcal{Z}}^{{\bf x},\varepsilon}_t  = 0\}, & 
{\rm if }~ \inf\{t \geq 0: {\mathcal{Z}}^{{\bf x},\varepsilon}_t  = 0\} < \inf\{t \geq 0: {\mathcal{Z}}^{{\bf  x},\varepsilon}_t  = r\}, \\ \infty, & {\rm otherwise}. \end{cases}
\]
Since $M^\varepsilon$ is a small perturbation of $M$ on $S \times [0,r]$
(formula (\ref{genpp0})), from (\ref{bnn}) it follows that, for all sufficiently small $\varepsilon$,
\begin{equation} \label{bnn2}
\mathrm{P}\left({\eta}^{{\bf x}, \varepsilon}_r < T,~\|\mathcal{X}^{{\bf x}}_{{\eta}^{{\bf x}}} - 
\mathcal{X}^{{\bf x},\varepsilon}_{{\eta}^{{\bf x},\varepsilon}_r} \| \leq \eta \right) \geq 1-  2\eta,~~~{\bf x} \in B_n.
\end{equation}
Observe that
\[
|u^\varepsilon({\bf x}) - \mathrm{E} \left( f(\mathcal{X}^{{\bf x},\varepsilon}_{{\eta}^{{\bf x},\varepsilon}_r});~ 
{\eta}^{{\bf x}, \varepsilon}_r < T \right) | \leq 2 \eta \sup |f|,
\]
\[
|u({\bf x}) - \mathrm{E} \left(
f(\mathcal{X}^{{\bf x}}_{{\eta}^{{\bf x}}}),~{\eta}^{{\bf x}, \varepsilon}_r < T \right) |\leq 2 \eta \sup |f|. 
\]
Therefore, from (\ref{bnn2}) it follows that, for all sufficiently small $\varepsilon$,
\[
 |u^\varepsilon({\bf x}) -  u({\bf x}) | \leq \mathrm{E} \left( |f(\mathcal{X}^{{\bf x},\varepsilon}_{{\eta}^{{\bf x},\varepsilon}_r}) - 
f(\mathcal{X}^{{\bf x}}_{{\eta}^{{\bf x}}}) |;~ 
{\eta}^{{\bf x}, \varepsilon}_r < T \right) | + 4 \eta \sup |f| \leq
\]
\[
\sup_{x_1, x_2 \in S, \|x_1 - x_2\| \leq \eta} |f(x_1)  -f(x_2)| +  8 \eta \sup |f| \leq \frac{\delta}{2},~~~{\bf x} \in B_n,
\]
where the last inequality is obtained by selecting a sufficiently small $\eta$. Thus,
\[
|u^\varepsilon({\bf x}) -  \overline{u}| \leq | u({\bf x}) - \overline{u}| +  |u^\varepsilon({\bf x}) -  u({\bf x}) |  \leq \delta, ~~~{\bf x} \in B_n.
\]

Next, consider the case when $S$ is repelling. Again, let $n \in \mathbb{N}$ be such that $|u({\bf x}) - \overline{u}| \leq \delta/2$ for
${\bf x} \in B_n$. Fix $\eta > 0$, to be specified later.  For $r > 0$, let
\[
{\eta}^{{\bf x}}_r =  \begin{cases} \inf\{t \geq 0: {\mathcal{Z}}^{{\bf x}}_t  = 0\}, & 
{\rm if }~ \inf\{t \geq 0: {\mathcal{Z}}^{{\bf x}}_t  = 0\} < \inf\{t \geq 0: {\mathcal{Z}}^{{\bf x}}_t  = r\}, \\ \infty, & {\rm otherwise}. \end{cases}
\]
Let $r$ be sufficiently large so that 
\[
|u({\bf x}) - \mathrm{E} \left( f(\mathcal{X}^{{\bf x}}_{{\eta}^{{\bf x}}_r}) | {\eta}^{{\bf x}}_r < \infty \right)| \leq \eta
\]
for all ${\bf x} \in B_n$. From the proximity of $\mathcal{X}^{{\bf x}}_t$ and $\mathcal{X}^{{\bf x},\varepsilon}_t$, 
using the same arguments as above, it is easy to show that 
\[
|\mathrm{E} \left( f(\mathcal{X}^{{\bf x},\varepsilon}_{{\eta}^{{\bf x},\varepsilon}_r}) | {\eta}^{{\bf x},\varepsilon}_r < \infty \right) -  
\mathrm{E} \left( f(\mathcal{X}^{{\bf x}}_{{\eta}^{{\bf x}}_r}) | {\eta}^{{\bf x}}_r < \infty \right) | \leq \eta,~~~{\bf x} \in B_n, 
\]
for all sufficiently small $\varepsilon$, and therefore 
$|\mathrm{E} \left( f(\mathcal{X}^{{\bf x},\varepsilon}_{{\eta}^{{\bf x},\varepsilon}_r}) | {\eta}^{{\bf x},\varepsilon}_r < \infty \right)  - u({\bf x})| \leq 2 \eta$.  Representing  
$u^\varepsilon(x) =  \mathrm{E} f(X^{x,\varepsilon}_{ \tau^{x,\varepsilon} (S)})$ as a sum of contributions
from successive visits to $B_n$ after reaching the surface defined by $\{{\bf z} = r\}$, we obtain 
$|u^\varepsilon({\bf x})   - u({\bf x})| \leq 2 \eta \leq \delta/2$, where the last inequality follows by taking a sufficiently small $\eta$. Thus,
$|u^\varepsilon({\bf x}) -  \overline{u}| \leq \delta$ for ${\bf x} \in B_n$. \qed


\subsection{The measure induced by the process stopped at $\Gamma_\varkappa$}

Let ${\tilde{\nu}}^{x,\varepsilon}_{\varkappa}$ 
be the measure on $\Gamma_\varkappa$ induced by 
$X^{x,\varepsilon}_{ \tau^{x,\varepsilon} (\Gamma_\varkappa)}$, where $x \in S$.   Since there is a natural bijection between 
$\Gamma_\varkappa$ and $S$ for each sufficiently small $\varkappa$, $\tilde{\nu}^{x,\varepsilon}_\varkappa$ can also be viewed as a measure on $S$. We are interested in the asymptotics  of $\tilde{\nu}^{x,\varepsilon}_\varkappa$ as $\varepsilon \downarrow 0$ first and $\varkappa \downarrow 0$ next. The arguments leading to the following result are very similar to those in Section~\ref{sseca}, so we do not provide a proof here.  
\begin{lemma} \label{indmo} For each sufficiently small $\varkappa > 0$, there is a measure $\tilde{\nu}_\varkappa$ such that, for each
$f \in C(\Gamma_\varkappa)$,
\[
\lim_{\varepsilon \downarrow 0} \int_{\Gamma_\varkappa} f d \tilde{\nu}^{x,\varepsilon}_\varkappa  = \int_{\Gamma_\varkappa} f d \tilde{\nu}_\varkappa
\]
uniformly in $x \in S$. There is a measure $\tilde{\nu}$ on $S$ such that $\tilde{\nu}_\varkappa \Rightarrow \tilde{\nu}$ as $\varkappa \downarrow 0$.
\end{lemma}
The limiting measure $\tilde{\nu}$ can be identified as follows.  Let $\tilde{X}^x_t = (\tilde{Y}^x_t, \tilde{Z}^x_t)$ be the
family of diffusion processes on $S \times (0, \infty)$ with the generator 
\[
K u  = L_y u  +  z^2 \alpha(y) \frac{\partial^2 u}{\partial z^2} +   z  \beta(y) \frac{\partial u}{\partial z}  + 
z {\mathcal{D}}_y \frac{\partial u}{\partial z}
\]
defined in Remark 2 after Lemma~\ref{gene}.
This operator already came up in the proof of Lemma~\ref{spec}; it can be obtained from the generator of  ${X}^x_t$ by discarding the last term in~(\ref{opr}) (this term is discarded since its contribution becomes insignificant near $S$). 

As in the proof of Lemma~\ref{sttm}, let $\Phi(x) = \Phi(y,z) = (y,\ln({z}))$ be the mapping from $S \times (0,\infty)$ to 
$S \times\mathbb{R}$.
The generator of the process $ \Phi(\tilde{X}^x_t) =   (\tilde{Y}^x_t, \ln(\tilde{Z}^x_t))$ 
on $S \times \mathbb{R}$ is the operator 
\[
\tilde{\mathcal{A}} u   = L_y u  +  \alpha(y) (\frac{\partial^2 u}{\partial  { z}^2} - 
\frac{\partial u}{\partial { z}} ) +
\beta(y) \frac{\partial u}{\partial  { z}}  + 
{\mathcal{D}}_y \frac{\partial u}{\partial  { z}}
\]
(compare with (\ref{genaa})).
In the case when $\overline{\alpha} < \overline{\beta}$, 
let $\hat{\nu}^x_\varkappa$ be the measure on $\Phi(\Gamma_\varkappa)$ induced by the process $\Phi(\tilde{X}^x_t)$ 
stopped when it hits $\Phi(\Gamma_\varkappa)$. In the case when $\overline{\alpha} > \overline{\beta}$, 
$\hat{\nu}^x_\varkappa$ is defined in the same 
way, but with the process conditioned on  reaching $\Phi(\Gamma_\varkappa)$.

It turns out that, for each $\varkappa > 0$, 
$\hat{\nu}^x_\varkappa \Rightarrow \tilde{\nu}$ as $x \rightarrow S$, i.e., when the initial point $(y,\ln(z))$ 
of the process $\Phi(\tilde{X}^x_t)$  satisfies $\ln(z) \rightarrow -\infty$. 
The value of $\varkappa$ is not important here since the coefficients of $\tilde{\mathcal{A}} $ do not depend on $z$. Similarly to
Lemma~\ref{mleee},
for each $f \in C(S)$, the value of $\int_S f d \tilde{\nu}$ can be expressed in terms of solutions to PDEs with the operator $\tilde{\mathcal{A}}$ on the domain $\{(y,z): 
y \in S, z \leq -\ln(\varphi(y))/\gamma \}$. 

\subsection{The time it takes to leave a neighborhood of a repelling boundary and the time it takes to reach an attracting boundary}
\label{sserea}

%

We start with the case of a repelling boundary.
\begin{lemma} \label{lee1}   Suppose that $\gamma < 0$.  

(a) For each sufficiently small $\varkappa > 0$, there is ${\varepsilon_0 = \varepsilon_0(\varkappa) > 0}$ such that the random variables 
$\tau^{x,\varepsilon} (\Gamma_\varkappa)/ |\ln(\varepsilon)|$ are uniformly integrable in $x \in V_{0,\varkappa}$, 
$0 < \varepsilon \leq \varepsilon_0$.

(b) For each sufficiently small $\varkappa > 0$, there is ${\varepsilon_0 = \varepsilon_0(\varkappa) > 0}$ such that 
\[
\lim_{c \downarrow 0} \mathrm{P}\left( \frac{\tau^{x,\varepsilon} (\Gamma_\varkappa)}{ |\ln(\varepsilon)|} < c\right) = 0
\]  
uniformly in $x \in S$, 
$0 < \varepsilon \leq \varepsilon_0$. 

(c) For each sufficiently small $\varkappa > 0$ and $0 < \zeta < \varkappa$, 
there is ${\varepsilon_0 = \varepsilon_0(\varkappa, \zeta) > 0}$ such that
$\tau^{x,\varepsilon} ( \Gamma_{\varkappa} ) $
are uniformly integrable in $x \in \Gamma_{\zeta}$, $0 < \varepsilon \leq \varepsilon_0$.  

\end{lemma}
\proof For $n \in \mathbb{Z}$, define $G_n    =  
\{(y,z): \psi(y) + \ln(z) = n\}$ - these are the analogues of
the sets $B_n$ from (\ref{subsets}), but now in $(y,z)$ coordinates.  Also define $G_n^+    =  
\{(y,z): \psi(y) + \ln(z) \geq n\}$. For $x \in G_n$, let 
\[
{\sigma}^{{x},\varepsilon} = \inf\{t \geq 0: {{X}}^{ x,\varepsilon}_t \in G_{n-1} \bigcup G_{n+1}\}.
\]
Using the smallness of the coefficients of $R$ for small $z$ (see (\ref{opr})) and the fact that the relative contribution of 
$\varepsilon^2 \tilde{L}$ to the generator of ${{X}}^{ x,\varepsilon}_t $ is small away from the boundary, it is easy to show, similarly to the proof of (\ref{slnn}), that there are $c >0$,
$n_0 < 0$, and $r >0$ such that 
\begin{equation} \label{est1}
\mathrm{P} ({{X}}^{ x, \varepsilon}_{{\sigma}^{{x, \varepsilon}} }  \in G_{n+1}) \geq \frac{1}{2} + c,~~~
{ x} \in G_n,~~[\ln(r \varepsilon)] \leq n \leq n_0.
\end{equation}
Inductively define a sequence of stopping times ${\sigma}^{{x},\varepsilon}_k$ as follows: ${\sigma}^{{x},\varepsilon}_0 = 0$ and, assuming that ${{X}}^{ x, \varepsilon}_{{\sigma}^{{x, \varepsilon}}_k} \in G_n$, 
\[
{\sigma}^{{x},\varepsilon}_{k+1} = \inf\{t \geq {\sigma}^{{x},\varepsilon}_{k}: {{X}}^{ x,\varepsilon}_t \in G_{n-1} \bigcup G_{n+1}\}.
\]
Let $K = \min\{k: {{X}}^{ x, \varepsilon}_{{\sigma}^{{x, \varepsilon}}_k} 
\in G_{n_0+1}\}$.
 Then $({{X}}^{ x, \varepsilon}_{{\sigma}^{{x, \varepsilon}}_{k \wedge K}},  
{\sigma}^{{x, \varepsilon}}_{k \wedge K })$, $k \geq 0$, is a Markov renewal process on the state space
 $\bigcup_{n = -\infty}^{n_0+1}
G_n$. (We stop the process at $G_{n_0+1}$ since the sets $G_{n}$ are not defined
for large $n$).

The inter-arrival times (conditioned on starting at $x \in \bigcup_{ n = -\infty }^{n_0} G_n$) are distributed as 
${\sigma}^{{x},\varepsilon}$,  which are uniformly integrable in $x \in \bigcup_{ n = [\ln(r \varepsilon)] }^{n_0} G_n$,
$0 \leq \varepsilon \leq \varepsilon_0$ (with some $\varepsilon_0 >0$) and
satisfy  
\begin{equation} \label{trtxxx}
\lim_{c \downarrow 0} \mathrm{P}( {{\sigma}^{x,\varepsilon}  } < c) = 0
\end{equation}
uniformly in $x$, $\varepsilon$.
To see this, write the process $X^{x,\varepsilon}_t$ in the coordinate form as $(Y^{x,\varepsilon}_t, Z^{x,\varepsilon}_t)$ and consider
the process $(Y^{x,\varepsilon}_t, \ln(Z^{x,\varepsilon}_t))$. The generator of the latter process is uniformly elliptic and
its coefficients are bounded (uniformly in $[\ln(r \varepsilon)] \leq n \leq n_0$) 
 in the domain bounded by $\Phi(G_{n-1})$ and $\Phi(G_{n+1})$, where
$\Phi(y,z) = (y, \ln(z))$. This leads to the uniform integrability and (\ref{trtxxx}).

Fix a  small $\varkappa > 0$ such that $\Gamma_\varkappa$ is defined,
and choose $r >0$ and $n_0$ such that (\ref{trtxxx}) holds and $G_{n_0+1} \subset V_{0,\varkappa}$ ($V_{0,\varkappa}$ is defined in the
beginning of Section~\ref{neigh}). 
In time $\tau^{x,\varepsilon} (\Gamma_\varkappa)$ 
that it takes the process $X^{x,\varepsilon}_t$ starting at $x \in S$ to reach $\Gamma_\varkappa$, the corresponding Markov renewal process
makes at least $n_0 - [\ln(r \varepsilon)]$ transitions between different pairs of sets $G_{n}$ and 
$G_{n+1}$ with
$[\ln(r \varepsilon)] \leq n \leq n_0$. 
 The estimate  (\ref{trtxxx}) on the inter-arrival times now implies part (b) of the lemma.
We can also estimate the time $\tau^{x,\varepsilon} (\Gamma_\varkappa)$ from above. 
Indeed, $ \tau^{x,\varepsilon} (G_{[\ln(r \varepsilon)]}^+) $ is uniformly integrable in $x \in V_{0,\varkappa}$,  $0 < \varepsilon \leq \varepsilon_0$,  since the process $\mathcal{X}^{{\bf x},\varepsilon}_t$ defined in (\ref{chp}) is non-degenerate. Once the process ${X}^{{x},\varepsilon}_t$
reaches $G_{[\ln(r \varepsilon)]}^+$, there is a positive
probability that it reaches $ G_{n_0}$ (and, consequently, $\Gamma_\varkappa$ with sufficiently small $\varkappa$) in time that is logarithmic in $\varepsilon$ due to (\ref{est1}) and the uniform integrability of ${\sigma}^{{x},\varepsilon}$. 
Part (a) of the lemma now easily follows from the strong Markov property. Part (c) is similar to part (a) - the time it takes a Markov renewal
process with a positive drift to move a finite distance to the right is uniformly integrable if the transition times are uniformly integrable.
\qed
\\

Now let us consider the case of an attracting boundary. We can still define the Markov renewal process associated with the
diffusion $X^{x,\varepsilon}_t$. The inter-arrival times are still uniformly integrable and satisfy (\ref{trtxxx}),  but the drift is now directed
towards $S$, i.e., instead of (\ref{est1}), we have
\[
\mathrm{P} ({{X}}^{ x, \varepsilon}_{{\sigma}^{{x, \varepsilon}} }  \in G_{n-1}) \geq \frac{1}{2} + c,~~~
{ x} \in G_n,~~[\ln(r \varepsilon)] \leq n \leq n_0.
\]
Arguments similar to those used in the proof of Lemma~\ref{lee1} can be used to justify the following lemmas. 
\begin{lemma} \label{lee2}  Suppose that $\gamma > 0$.   

(a) For each sufficiently small $\varkappa > 0$, there is ${\varepsilon_0 = \varepsilon_0(\varkappa) > 0}$ such that the random variables 
$\tau^{x,\varepsilon} (S \bigcup \Gamma_\varkappa)/ |\ln(\varepsilon)|$ are uniformly integrable in $x \in V_{0,\varkappa}$, 
$0 < \varepsilon \leq \varepsilon_0$.

(b) For each sufficiently small $\varkappa > 0$, there is ${\varepsilon_0 = \varepsilon_0(\varkappa) > 0}$ such that 
\[
\lim_{c \downarrow 0} \mathrm{P}\left( \frac{\tau^{x,\varepsilon} (S)}{ |\ln(\varepsilon)|} < c\right) = 0
\]  
uniformly in $x \in \Gamma_\varkappa$, 
$0 < \varepsilon \leq \varepsilon_0$. 

(c) For each $r > 0$ and each sufficiently small $\varkappa$, there is ${\varepsilon_0 = \varepsilon_0(\varkappa, r) > 0}$ such that
$\tau^{x,\varepsilon} ( S \bigcup  \Gamma_{\varkappa} ) $
are uniformly integrable in $x \in \Gamma_{r \varepsilon}$, $0 < \varepsilon \leq \varepsilon_0$.  
\end{lemma}
The next lemma implies that the process starting in a small neighborhood of $S$ is likely to reach $S$ before leaving a much larger neighborhood of $S$. 
\begin{lemma} \label{trp5}
 Suppose that $\gamma > 0$.  

(a) For each sufficiently small $\varkappa > 0$, 
\[
\lim_{\zeta, \varepsilon \downarrow 0} \mathrm{P} (X^{x, \varepsilon}_{ \tau^{x,\varepsilon} ( S \bigcup  \Gamma_{\varkappa} )} 
\in S) = 1
\]
uniformly in $x \in V_{0,\zeta}$.

(b) For each $r > 0$, 
\[
\lim_{s \uparrow \infty, \varepsilon \downarrow 0} \mathrm{P} (X^{x, \varepsilon}_{ \tau^{x,\varepsilon} ( S \bigcup  \Gamma_{s \varepsilon} )} 
\in S) = 1
\]
uniformly in $x \in V_{0,r \varepsilon}$.
\end{lemma}

\subsection{The time it takes to leave a neighborhood of an attracting boundary and the time it takes to reach a repelling boundary}

This time, we start with the case of an attracting boundary. Let us estimate the probability that
the process starting near $S$ escapes from a small neighborhood of $S$ without visiting $S$.
\begin{lemma} \label{prat} Suppose that $\gamma > 0$. For each $\eta > 0$,  for all sufficiently small $\varkappa$ and sufficiently 
large $s > 0$  (both dependent on $\eta$), there is $\varepsilon_0 > 0$ such that
\begin{equation} \label{prex2aa}
(1 - \eta)  \frac{ \zeta^\gamma}{\varkappa^\gamma} \leq \mathrm{P} 
( X^{x,\varepsilon}_{ \tau^{x,\varepsilon} (S \bigcup \Gamma_\varkappa)} \in \Gamma_{\varkappa}) \leq (1 + \eta)
  \frac{ \zeta^\gamma}{\varkappa^\gamma},
\end{equation} 
provided that $x \in \Gamma_\zeta$ with $s \varepsilon \leq \zeta \leq \varkappa$ and $0 < \varepsilon \leq \varepsilon_0$.
\end{lemma}
\proof
First, we will consider the behavior of the process in $V_{r \varepsilon, \varkappa}$ with $r, \varkappa > 0$. We claim that, given  $\eta > 0$, for all sufficiently large $r$ and sufficiently small $\varkappa$, 
\begin{equation} \label{prex}
(1 - \eta)  \frac{\zeta^\gamma - (r \varepsilon)^\gamma}{\varkappa^\gamma} \leq \mathrm{P} 
( X^{x,\varepsilon}_{ \tau^{x,\varepsilon} (\Gamma_{r \varepsilon} \bigcup \Gamma_\varkappa)} \in \Gamma_{\varkappa}) \leq (1 + \eta)
\frac{\zeta^\gamma - (r \varepsilon)^\gamma}{\varkappa^\gamma},
\end{equation}  
provided that $x \in \Gamma_{\zeta}$ with $2 r \varepsilon \leq \zeta \leq \varkappa$ and  $\varepsilon$ is sufficiently small. The main idea is that the region $V_{r \varepsilon, \varkappa}$
is separated from $S$ (so that the perturbation $\varepsilon^2   \tilde{L}$ to the operator $L$ does not play a big role) and, at the same time, is sufficiently close to $S$ so that only the leading terms of the expansion of the coefficients of $L$ play a role (i.e., the operator $R$ in (\ref{opr}) is sufficiently small). Thus 
the process  $\varphi(Y^{x, \varepsilon}_t)(Z^{x, \varepsilon}_t)^\gamma$ is nearly a martingale in $V_{r \varepsilon, \varkappa}$ (or, equivalently,
the generator of the process applied to $\varphi(y)z^\gamma$ is nearly zero), which suggests that (\ref{prex}) should hold.

Let us now give a rigorous argument. Recall from the proof of Lemma~\ref{spec} that the top eigenvalue of the operator 
in the left hand side of (\ref{mg}) satisfies $\lambda_0 = 0$, $\lambda'_0 < 0$, and $\lambda_\gamma = 0$. Since the top
eigenvalue depends continuously on the parameter,  by the uniqueness part of Lemma~\ref{spec}, the eigenvalues corresponding to the values of the parameter that are slightly 
smaller than $\gamma$ are negative, and thus
 there exist $\gamma_1 \in (\max(0, \gamma-1) ,\gamma)$ and a positive-valued function 
$\varphi_1  \in C^1(S)$ satisfying $\int_S \varphi_1 d \pi = 1$ such that
\begin{equation} \label{yyu1}
L_y \varphi_1 + \alpha \gamma_1 (\gamma_1 -1) \varphi_1 + \beta \gamma_1 \varphi_1 + \gamma_1 \mathcal{D}_y \varphi_1 = -c_1 \varphi_1,
\end{equation}
where $c_1 >0$.  Similarly, there exist  $\gamma_2 \in (\gamma, \gamma +1)$ and a positive-valued function 
$\varphi_2  \in C^1(S)$ satisfying $\int_S \varphi_2 d \pi = 1$ such that
\begin{equation} \label{yyu2}
L_y \varphi_2 + \alpha \gamma_2 (\gamma_2 -1) \varphi_2 + \beta \gamma_2 \varphi_2 + \gamma_2 \mathcal{D}_y \varphi_2 = c_2 \varphi_2,
\end{equation}
where $c_2 >0$. Let
\[
u(y,z) = \varphi(y) z^\gamma,~~ u_1(y,z)  = \varphi_1(y) z^{\gamma_1},~~u_2(y,z)  = \varphi_2(y) z^{\gamma_2}.
\]
By Lemma~\ref{gene} and (\ref{mg}), the function 
\begin{equation} \label{fnff}
f^\varepsilon(y,z) = (u(y,z) - (r\varepsilon)^\gamma)/(\varkappa^\gamma - (r\varepsilon)^\gamma)
\end{equation}
 satisfies
\[
L^\varepsilon f^\varepsilon (y,z) = (R +\varepsilon^2   \tilde{L})f^\varepsilon (y,z),~~~(y,z) \in V_{r \varepsilon, \varkappa},
\]
\[
f^\varepsilon|_{\Gamma_{r \varepsilon}} = 0,~~~f^\varepsilon|_{\Gamma_{\varkappa}} = 1.
\]
The function $v^\varepsilon(x) = \mathrm{P} 
( X^{x,\varepsilon}_{ \tau^{x,\varepsilon} (\Gamma_{r \varepsilon} \bigcup \Gamma_\varkappa)} \in \Gamma_{\varkappa})$, which is what we are interested in, satisfies the same boundary conditions and
the same equation but with zero instead of $(R +\varepsilon^2   \tilde{L})f^\varepsilon (y,z)$ in the right hand side. Thus, 
\begin{equation} \label{probn}
\mathrm{P} 
( X^{x,\varepsilon}_{ \tau^{x,\varepsilon} (\Gamma_{r \varepsilon} \bigcup \Gamma_\varkappa)} \in \Gamma_{\varkappa}) = v^\varepsilon(x) =  
f^\varepsilon(x) - g^\varepsilon(x),
\end{equation}
 where $g^\varepsilon$ solves
\begin{equation} \label{sl11}
L^\varepsilon g^\varepsilon (y,z) = (R +\varepsilon^2   \tilde{L})f^\varepsilon (y,z),~~~(y,z) \in V_{r \varepsilon, \varkappa},
\end{equation}
\[
g^\varepsilon|_{\Gamma_{r \varepsilon}} = 0,~~~g^\varepsilon|_{\Gamma_{\varkappa}} = 0.
\]
Observe that there is $C > 0$ such that 
\begin{equation} \label{rs1s}
|(R +\varepsilon^2   \tilde{L})f^\varepsilon (y,z)| \leq  C \varkappa^{-\gamma} (z^{\gamma+1} + \varepsilon^2 z^{\gamma-2})
\end{equation}
for all sufficiently small $\varepsilon$. By Lemma~\ref{gene}, using (\ref{yyu1}), we can find an arbitrarily small $k_1 > 0$ (by taking $r$ sufficiently large and $\varkappa$
sufficiently small) such that
\[
k_1 \varepsilon^{ \gamma  - \gamma_1 } L^\varepsilon ( u_1(y,z)) \leq - C   \varepsilon^2 z^{ \gamma-2},~~~(y,z) \in V_{r \varepsilon, \varkappa},
\] 
for all sufficiently small $\varepsilon$. Similarly,  using (\ref{yyu2}), we can find an arbitrarily small $k_2 > 0$  (by taking $r$ sufficiently large and $\varkappa$
sufficiently small) such that
\[
 k_2  L^\varepsilon ( u_2(y,z)) \geq  C    z^{\gamma+1},~~~(y,z) \in V_{r \varepsilon, \varkappa},
\] 
for all sufficiently small $\varepsilon$. By taking $k_2$ sufficiently small,  we can find an arbitrarily small $k >0$ such that
\[
k  \inf_{(y,z) \in {\Gamma_{r\varepsilon}} } u(y,z) \geq k_2  \sup _{(y,z) \in {\Gamma_{r\varepsilon}} } u_2(y,z),~~~
k  \inf_{(y,z) \in {\Gamma_\varkappa} } u(y,z) \geq k_2  \sup _{(y,z) \in {\Gamma_\varkappa} } u_2(y,z)
\]
for all sufficiently small $\varepsilon$. Thus the function $\tilde{g}^\varepsilon  = \varkappa^{-\gamma}(k_1 \varepsilon^{
\gamma  - \gamma_1} u_1 + k u -
k_2 u_2)$ satisfies
\[
L^\varepsilon \tilde{g}^\varepsilon (y,z) \leq  - C (1-k) \varkappa^{-\gamma} (z^{\gamma+1} + \varepsilon^2 z^{ \gamma-2 }),~~~(y,z) \in V_{r \varepsilon, \varkappa},
\]
\[
\tilde{g}^\varepsilon|_{\Gamma_{r \varepsilon}} \geq 0,~~~\tilde{g}^\varepsilon|_{\Gamma_{\varkappa}} \geq 0.
\]
Comparing this with (\ref{sl11}), (\ref{rs1s}) and using the stochastic representation for the solutions $\tilde{g}^\varepsilon$ and ${g}^\varepsilon$
of the respective equations, we obtain that 
\[
\tilde{g}^\varepsilon(y,z) \geq (1-k) |{g}^\varepsilon(y,z)|, ~~~(y,z) \in V_{r \varepsilon, \varkappa}.
\]
Thus
\[
 |{g}^\varepsilon(y,z)| \leq \frac{\varkappa^{-\gamma}}{1 - k} (k_1 \varepsilon^{\gamma  - \gamma_1} \varphi_1(y) z^{\gamma_1} + k \varphi(y) z^{\gamma} - k_2 \varphi_2(y) z^{\gamma_2})  
\]
\[
\leq \frac{\varkappa^{-\gamma}}{1 - k} (k_1 \varepsilon^{\gamma  - \gamma_1} \varphi_1(y) z^{\gamma_1} + k \varphi(y) z^{\gamma} ).
\]
For each $\tilde{\eta} > 0$, by making $k_1$ and $k$ sufficiently small, we can make sure that
\[
\sup_{(y,z) \in \Gamma_{\zeta}} |{g}^\varepsilon(y,z)| \leq   \frac{\tilde{\eta} \zeta^\gamma}{\varkappa^\gamma} 
\]
for all sufficiently small $\varepsilon$. Combining this with (\ref{probn}) and the definition (\ref{fnff}) of $f^\varepsilon$,
we obtain
\[
 \frac{\zeta^\gamma - (r \varepsilon)^\gamma}{\varkappa^\gamma - (r \varepsilon)^\gamma} - \frac{\tilde{\eta} \zeta^\gamma}{\varkappa^\gamma} \leq \mathrm{P} 
( X^{x,\varepsilon}_{ \tau^{x,\varepsilon} (\Gamma_{r \varepsilon} \bigcup \Gamma_\varkappa)} \in \Gamma_{\varkappa}) \leq 
\frac{\zeta^\gamma - (r \varepsilon)^\gamma}{\varkappa^\gamma- (r \varepsilon)^\gamma} + \frac{\tilde{\eta} \zeta^\gamma}{\varkappa^\gamma}.
\]
Since $\tilde{\eta}$ can be taken arbitrarily small and 
 $\zeta \geq 2r \varepsilon$,  this estimate implies that (\ref{prex}) holds for all sufficiently small $\varepsilon$.

In order to prove the first inequality in (\ref{prex2aa}), we use the first inequality in
(\ref{prex}), which implies that
\[
(1 - \eta)  \frac{\zeta^\gamma - (r \varepsilon)^\gamma}{\varkappa^\gamma}\leq
\mathrm{P} 
( X^{x,\varepsilon}_{ \tau^{x,\varepsilon} (S \bigcup \Gamma_\varkappa)} \in \Gamma_{\varkappa}) 
\]
for $x \in \Gamma_{\zeta}$. It remains to note that $(\zeta^\gamma - (r \varepsilon)^\gamma)/\zeta^\gamma$  can
be made arbitrarily close to one for $\zeta \geq s \varepsilon$ by selecting a sufficiently large $s$. 

Next, let us prove the second inequality in  (\ref{prex2aa}). The process $X^{x,\varepsilon}_t$ may
reach $\Gamma_\varkappa$ either before visiting $\Gamma_{r \varepsilon}$ or after visiting $\Gamma_{r \varepsilon}$ and returning to $\Gamma_{s \varepsilon}$.
Therefore, by the strong Markov property of the process, 
\[
\sup_{x \in \Gamma_{\zeta}} \mathrm{P} 
( X^{x,\varepsilon}_{ \tau^{x,\varepsilon} (S \bigcup \Gamma_\varkappa)} \in \Gamma_{\varkappa}) 
\]
\begin{equation} \label{inq1}
\leq  \sup_{x \in \Gamma_{\zeta}} \mathrm{P} 
( X^{x,\varepsilon}_{ \tau^{x,\varepsilon} (\Gamma_{r \varepsilon} \bigcup \Gamma_\varkappa)} \in \Gamma_{\varkappa}) 
\end{equation}
\[
+ \sup_{x \in \Gamma_{r  \varepsilon}}  \mathrm{P} ( X^{x,\varepsilon}_{ \tau^{x,\varepsilon} (S \bigcup \Gamma_{\zeta})} \in \Gamma_{\zeta})
\sup_{x \in \Gamma_{\zeta}} \mathrm{P} 
( X^{x,\varepsilon}_{ \tau^{x,\varepsilon} (S \bigcup \Gamma_\varkappa)} \in \Gamma_{\varkappa}). 
\]
By part (b) of Lemma~\ref{trp5},
given  $\eta > 0$, for all sufficiently large $s$,
\begin{equation} \label{inq2}
\sup_{x \in \Gamma_{r  \varepsilon}}  \mathrm{P} ( X^{x,\varepsilon}_{ \tau^{x,\varepsilon} (S \bigcup \Gamma_{\zeta})} \in \Gamma_{\zeta}) \leq  \eta, 
\end{equation}
provided that $\zeta \geq s \varepsilon$ and $\varepsilon$ is  sufficiently small. Combining (\ref{inq1}) and (\ref{inq2}) with the second inequality in 
(\ref{prex}), we obtain the second estimate in  in  (\ref{prex2aa}).
\qed
\\

Now we are ready to prove the result concerning the exit time from a neighborhood of an attracting boundary.
\begin{lemma} \label{mltime} Suppose that $\gamma > 0$.  

(a) There is a positive constant $q$ such that
for each $\eta > 0$  and each sufficiently small $\varkappa$ (depending on $\eta$), there is $\varepsilon_0 = \varepsilon_0(\varkappa) > 0 $ such that
\begin{equation} \label{time1}
 q (1 - \eta)   (\frac{\varkappa}{\varepsilon})^\gamma \leq \mathrm{E} \tau^{x,\varepsilon} (\Gamma_\varkappa) \leq q (1 + \eta)   (\frac{\varkappa}{\varepsilon})^\gamma,
\end{equation}
for all $x \in S$, $0 < \varepsilon \leq \varepsilon_0$.

(b)  For each sufficiently small $\varkappa$ (depending on $\eta$), there is $\varepsilon_0 = \varepsilon_0(\varkappa) > 0 $ such that  the random variables 
$\varepsilon^\gamma \tau^{x,\varepsilon} (\Gamma_\varkappa) $ are uniformly integrable in $x \in S$, 
$0 < \varepsilon \leq \varepsilon_0$.

(c) For each sufficiently small $\varkappa$ (depending on $\eta$), there is $\varepsilon_0 = \varepsilon_0(\varkappa) > 0 $ such that
\begin{equation} \label{bfb}
\lim_{c \downarrow 0} \mathrm{P}\left( \varepsilon^\gamma   {\tau^{x,\varepsilon} (\Gamma_\varkappa)}  < c\right) = 0
\end{equation}
uniformly in $x \in S$, 
$0 < \varepsilon \leq \varepsilon_0$. 
\end{lemma} 
\proof 
Consider the Markov renewal process $(\xi^{x, \varepsilon}_n, \sigma^{x, \varepsilon}_n)$ on the state space 
$M_\varkappa = S \bigcup \Gamma_\varkappa$ with the starting point $x \in M_\varkappa$. The process, which depends on the parameter $s > 0$, to be selected later, 
is defined as follows. For $n = 0$,  $\xi^{x, \varepsilon}_0 = x$, 
$\sigma^{x, \varepsilon}_0 = 0$. For $n \geq 1$, let
\[
\tilde{\sigma}^{x, \varepsilon}_n = \inf\{t > \sigma^{x, \varepsilon}_{n-1}: X^{x, \varepsilon}_t \in \Gamma_{s \varepsilon} \}.
\] 
Then we define
\[
{\sigma}^{x, \varepsilon}_n = \inf\{t > \tilde{\sigma}^{x, \varepsilon}_{n}: X^{x, \varepsilon}_t \in M_\varkappa \},~~~\xi^{x, \varepsilon}_n = 
X^{x,\varepsilon}_{{\sigma}^{x, \varepsilon}_n}.
\] 
Denote the transition function of $\xi^{x, \varepsilon}_n$ by $Q^\varepsilon(x, A)$, $x \in M_\varkappa$, $A \in \mathcal{B}(M_\varkappa)$. Let $\mu^{s,\varepsilon}$ be
the invariant measure for the Markov chain $\xi^{x, \varepsilon}_n$. 
Arguing as in the proof of Lemma~\ref{mlee1}, it is easy to see that, for each $s > 0$,
$\mu^{s,\varepsilon}$ converges weakly, as $\varepsilon \downarrow 0$, to some measure $\mu^s$ concentrated on $S$.

Define the random variable $N^{x,\varepsilon} $ as
\[
N^{x,\varepsilon} = \min\{n: \xi^{x, \varepsilon}_n \in \Gamma_\varkappa \}.
\]
Thus $\tau^{x,\varepsilon} (\Gamma_\varkappa) = {\sigma}^{x, \varepsilon}_{N^{x,\varepsilon} }$.  We make the following observations.

(a) Due to (\ref{prex2aa}),   $\varepsilon^{-\gamma} \mathrm{P} 
(\xi^{x, \varepsilon}_1 \in \Gamma_\varkappa)$ is bounded between two positive constants uniformly in $x \in S$ and $\varepsilon$, for each 
fixed value of $s$. 
Therefore, $\lim_{\varepsilon \downarrow 0} N^{x,\varepsilon} = \infty$ in probability (the process
$X^{x, \varepsilon}_t$ makes a large number excursions from $S$ to $\Gamma_{s \varepsilon}$ and back before reaching $ \Gamma_\varkappa$).

(b) Due to non-degeneracy of the process $\mathcal{{X}}^{\bf x,\varepsilon}_t$ (defined in the beginning of Section~\ref{sseca}) in 
a vicinity of $S$, the distribution of $\xi^{x, \varepsilon}_1$ has a density on $S$ that is bounded from above and below by positive
constants that are independent of $x \in S$ and $\varepsilon$. 

(c) The random variables $\tau^{x,\varepsilon} (\Gamma_{s \varepsilon})$, $x \in S$, 
$\varepsilon >0$, and 
$\tau^{x,\varepsilon} (S \bigcup \Gamma_\varkappa)$, $x \in \Gamma_{s \varepsilon}$, $\varepsilon >0$, 
are uniformly integrable (due to non-degeneracy of $\mathcal{{X}}^{\bf x,\varepsilon}_t$ and by Lemma~\ref{lee2}, part (c), respectively).
Therefore, ${\sigma}^{x, \varepsilon}_1$ are uniformly integrable in $x \in S$ and $\varepsilon > 0$.

From (a)-(c) it follows that the expected time  to reach $\Gamma_\varkappa$ is asymptotically equivalent to the expected
number of excursions between $M_\varkappa$ and $S_{s \varepsilon}$ multiplied by the expectation of the time spent on one excursion, where
the latter expectations are taken using the invariant measure on $M_\varkappa$ as the initial distribution. Note that the number of
excursions is, asymptotically, a geometric random variable with a parameter that tends to zero, and the expected number of excursions is 
the inverse of the probability that the process reaches $\Gamma_\varkappa$ during one excursion. Thus
\begin{equation} \label{jik}
\mathrm{E} \tau^{x,\varepsilon} (\Gamma_\varkappa) = \mathrm{E} {\sigma}^{x, \varepsilon}_{N^{x,\varepsilon} } \sim 
\frac{\int_S \mathrm{E} {\sigma}^{{\rm x}, \varepsilon}_{1} d \mu^{s,\varepsilon}({ \rm x})}{ \int_S \mathrm{P} 
(\xi^{{\rm x}, \varepsilon}_1 \in \Gamma_\varkappa)   d \mu^{s,\varepsilon}({\rm x})},~~~{\rm as}~~\varepsilon \downarrow 0,
\end{equation}
uniformly in $x \in S$. Here, we used symbol ${\rm x}$ in the right hand side to distinguish the variable used in the integral from the initial point $x \in S$. 

Observe that $\lim_{\varepsilon \downarrow 0} \mathrm{E} \tau^{x,\varepsilon} (\Gamma_{s \varepsilon}) = f(s,x)$ uniformly in $x \in S$ for some
positive function~$f$.  Moreover, $f$ is continuous in $x$ for each $s$ due to non-degeneracy of the process $\mathcal{{X}}^{\bf x,\varepsilon}_t$ defined in the beginning of Section~\ref{sseca}.

Applying arguments similar to those in Lemma~\ref{lee1}, it is easy to see that, due to the
boundary being attracting, it takes much longer
for the process starting at $S$ to reach $\Gamma_{s \varepsilon}$ than for the process starting at $\Gamma_{s \varepsilon}$ to reach $S 
\bigcup \Gamma_\varkappa$,
provided that $s$ is sufficiently large. More precisely, for each $\eta_1 > 0$, for all sufficiently large $s$ and sufficiently small $\varepsilon$,
\[
\sup_{x \in \Gamma_{s \varepsilon}} \mathrm{E} \tau^{x,\varepsilon} (S \bigcup \Gamma_\varkappa) \leq \eta_1 \int_S f(s,{\rm x}) d 
\mu^{s,\varepsilon}({\rm x}).
\]

Therefore, by the strong Markov property, given $\eta > 0$, for sufficiently small $\varkappa$ and sufficiently large $s$, 
\[ 
\int_S f(s,{\rm x}) d \mu^{s,\varepsilon}({\rm x}) \leq \int_S \mathrm{E} {\sigma}^{{\rm x}, \varepsilon}_{1} d \mu^{s,\varepsilon}( {\rm x})
\leq   (1+ \eta_1) \int_S 
f(s,{\rm x}) d \mu^{s,\varepsilon}({\rm x}),
\]
provided that $\varepsilon$ is sufficiently small, which implies that
\begin{equation} \label{smep}
(1 -\eta_1)\int_S f(s,{\rm x}) d \mu^{s }({\rm x}) \leq \int_S \mathrm{E} {\sigma}^{{\rm x}, \varepsilon}_{1} d \mu^{s,\varepsilon}( {\rm x})
\leq   (1+ 2\eta_1) \int_S 
f(s,{\rm x}) d \mu^{s }({\rm x}),
\end{equation}
provided that $\varepsilon$ is sufficiently small.
  For the denominator in the right-hand side of (\ref{jik}), by the strong Markov property and (\ref{prex2aa}), for an arbitrarily small $\eta_2 > 0$, we can
write 
\begin{equation} \label{ui2}
(1 - \eta_2)  \frac{ s^\gamma \varepsilon^\gamma}{\varkappa^\gamma} \leq \int_S \mathrm{P} 
(\xi^{{\rm x}, \varepsilon}_1 \in \Gamma_\varkappa)   d \mu^{s,\varepsilon}({\rm x}) \leq (1 + \eta_2)
  \frac{ s^\gamma \varepsilon^\gamma}{\varkappa^\gamma}.
\end{equation}
Combining (\ref{smep}) with (\ref{ui2}) and using (\ref{jik}), we obtain (\ref{time1}) with $q = \int_S 
f(s,{\rm x}) d \mu^{s }({\rm x})/s^\gamma$. It turns out that
$q$ does not depend on $s$ since $\mathrm{E} \tau^{x,\varepsilon} (\Gamma_\varkappa)$ does not depend on $s$.

Consider fixed $\varkappa$ and $s$, and let $\varepsilon \downarrow 0$. Since $\tau^{x,\varepsilon} (\Gamma_\varkappa) = {\sigma}^{x, \varepsilon}_{N^{x,\varepsilon} }$ is the first time when a Markov renewal process reaches the set $\Gamma_\varkappa$, and  $\lim_{\varepsilon \downarrow 0} Q^\varepsilon(x, \Gamma_\varkappa) = 0$,
uniformly in $x \in S$, from the uniform integrability of ${\sigma}^{x, \varepsilon}_1$  
and the boundedness from below of $\mathrm{E} {\sigma}^{x, \varepsilon}_1$, it follows that the random variables 
$\varepsilon^\gamma \tau^{x,\varepsilon} (\Gamma_\varkappa) $ are uniformly integrable in $x \in S$, 
$0 < \varepsilon \leq \varepsilon_0$, and (\ref{bfb}) holds.  \qed
\\

%
The next lemma shows that, at time scales $1 \ll t(\varepsilon) \ll  \varepsilon^{-\gamma}$, the distribution of 
$X^{x,\varepsilon}_{t(\varepsilon)}$ is close to the invariant measure $\pi$.  
\begin{lemma} \label{wcon} Suppose that $\gamma > 0$.  If $1 \ll t(\varepsilon) \ll  \varepsilon^{-\gamma}$ and $f \in C(\overline{D})$, 
then 
\[
\lim_{\varepsilon, \varkappa \downarrow 0} \sup_{x \in  V_{0,\varkappa}} |\mathrm{E} f (X^{x,\varepsilon}_{t(\varepsilon)}) - \int_S f d \pi| = 0.
\]
\end{lemma}
\proof Since the process $X^x_t$ is non-degenerate on $S$, it is exponentially mixing, and therefore, for each $\eta > 0$, one can find such $T > 0$ that
\[
\sup_{x \in  S} |\mathrm{E} f (X^{x}_{T}) - \int_S f d \pi| < \eta.
\]
Since $T$ is fixed, the same estimate, with $\eta$ replaced by $2 \eta$, holds for a small perturbation of the process.
Thus, there exist $\varepsilon_0 > 0$ and $\varkappa > 0$ such that 
\[
\sup_{0 < \varepsilon < \varepsilon_0, x \in  V_{0,\varkappa}} |\mathrm{E} f (X^{x,\varepsilon}_{T}) - \int_S f d \pi| < 2\eta.
\]
Lemmas~\ref{prat} and \ref{mltime} imply that, given $\eta > 0$, with probability larger than $1 - \eta$, it takes 
time of order $\varepsilon^{-\gamma}$ for the process to leave the small neighborhood $V_{0,\varkappa}$ of the attracting surface $S$, provided that it starts in a yet smaller
neighborhood $V_{0,\zeta}$ (with $\zeta$ that depends on $\varkappa$ and $\eta$). Therefore, considering the value of the
process at time $t(\varepsilon) - T$ as a new starting point and using the Markov property, we obtain
\[
\sup_{0 < \varepsilon < \varepsilon_0, x \in  V_{0,\zeta}} |\mathrm{E} f (X^{x,\varepsilon}_{t(\varepsilon)}) - \int_S f d \pi| < 3\eta,
\]
which implies the desired result.
\qed 
\\

Next, we discuss the case of a repelling boundary. We have the following counterpart of Lemma~\ref{prat}.
\begin{lemma} \label{prre} Suppose that $\gamma < 0$. There is $\rho > 0$ such that, for each $\eta > 0$,  for all sufficiently small $\varkappa$, sufficiently 
large $s_1 > 0$, and sufficiently small $s_2 > 0$  (all dependent of $\eta$), there is $\varepsilon_0 > 0$ such that
\begin{equation} \label{prex2aaz}
\rho (1 - \eta)  \frac{ \zeta^\gamma}{\varepsilon^\gamma} \leq \mathrm{P} 
( X^{x,\varepsilon}_{ \tau^{x,\varepsilon} (S \bigcup \Gamma_\varkappa)} \in S) \leq \rho (1 + \eta)
  \frac{ \zeta^\gamma}{\varepsilon^\gamma},
\end{equation} 
provided that $x \in \Gamma_\zeta$ with $s_1 \varepsilon \leq \zeta \leq s_2 \varkappa$ and $0 < \varepsilon \leq \varepsilon_0$.
\end{lemma}
Note the constant factor $\rho$ that depends on the coefficients of the perturbation $\tilde{L}$. Its presence is due to the fact that
the probability for the process $X^{x,\varepsilon}_t$ to reach $S$ 
prior to $\Gamma_\varkappa$ essentially depends on the behavior in an $\varepsilon$-neighborhood of $S$. The proof of Lemma~\ref{prre} is
similar to that of Lemma~\ref{prat}, so we do not provide it here.

\section{Asymptotics of transition probabilities between different components of the boundary} \label{atpr}
For $x \in \partial D$, inductively define a sequence of stopping times ${\sigma}^{{x},\varepsilon}_n$ as follows: 
${\sigma}^{{x},\varepsilon}_0 = 0$ and, assuming that ${{X}}^{ x, \varepsilon}_{{\sigma}^{{x, \varepsilon}}_n} \in S_k$, 
\[
{\sigma}^{{x},\varepsilon}_{n+1} = \inf\{t \geq {\sigma}^{{x},\varepsilon}_{n}: {{X}}^{ x,\varepsilon}_t \in \partial D \setminus S_k\}.
\]
 Then $({{X}}^{ x, \varepsilon}_{{\sigma}^{{x, \varepsilon}}_{n}},  
{\sigma}^{{x, \varepsilon}}_{n})$, $n \geq 0$, is a Markov renewal process on the state space
 $\partial D$. Let $Q^\varepsilon$ be its transition kernel. We are interested in the asymptotics, as $\varepsilon \downarrow 0$,
of $Q^\varepsilon(x, S_{j}) = \mathrm{P}(X^{ x, \varepsilon}_{{\sigma}^{x, \varepsilon}} \in S_{j})$ with 
$x \in S_{i}$, $i \neq j$, where we denote ${\sigma}^{x, \varepsilon} = {\sigma}^{x, \varepsilon}_1$. 
\begin{lemma} \label{trapl}
There exist constants $q_{i j} > 0$, $1 \leq i, j \leq m$, $i \neq j$, such that:

(a) If $S_{j}$ is attracting, then
\[
\lim_{\varepsilon \downarrow 0} Q^\varepsilon(x, S_{j}) = q_{i j},~~x \in S_{i}.
\]

(b) If $S_{j}$ is repelling and there is $k \neq i$ with $S_k$ attracting, then
\[
\lim_{\varepsilon \downarrow 0} \varepsilon^{\gamma_{j}} Q^\varepsilon(x, S_{j}) = q_{i j},~~x \in S_{i}.
\]

(c) If $S_k$ is repelling for each $k \neq i$, then take $k^*$ such that $\gamma_{k^*} \geq \max_{k \neq i} \gamma_k$.
In this case,
\[
\lim_{\varepsilon \downarrow 0} \left( \varepsilon^{\gamma_{j} - \gamma_{k^*}} Q^\varepsilon(x, S_{j}) \right) = q_{i j},~~x \in S_{i}.
\]
In each of the cases above, the convergence is uniform in $x \in S_{i}$. 
\end{lemma}
\proof The surfaces $\Gamma^k_\varkappa$, $1 \leq k \leq m$, are defined as in Section~\ref{neigh}, but with $S_k$ instead of a generic 
component of the boundary $S$.

(a)  Fix $\varkappa > 0$ sufficiently small so that the conclusion of Lemma~\ref{indmo} holds. 
Let $f, f^\varepsilon, g , g^\varepsilon \in C(\Gamma^{i}_\varkappa)$ be defined as
\[
f(x) = \mathrm{P}\left(\lim_{t \rightarrow \infty} {\rm dist}(X^x_t, S_{i})  = 0\right),~~~~
f^\varepsilon(x) = \mathrm{P}(X^{x,\varepsilon}_{ \tau^{x,\varepsilon}(\partial D)} \in S_{i}),
\]
\[
g(x) = \mathrm{P}\left(\lim_{t \rightarrow \infty} {\rm dist}(X^x_t, S_{j})  = 0\right),~~~~
g^\varepsilon(x) = \mathrm{P}(X^{x,\varepsilon}_{ \tau^{x,\varepsilon}(\partial D)} \in S_{j}) .
\]
Let $A^{x,\varepsilon}_n$, $n \geq 1$, $x \in S_{i}$, be the event that $X^{x, \varepsilon}_t$ makes at least $n$ transitions between
$S_{i}$ and $\Gamma^{i}_\varkappa$ prior to reaching $\partial D \setminus S_{i}$. Let $\tau^{x,\varepsilon}_n$ be the earliest time
when $X^{x, \varepsilon}_t$ visits $\Gamma^{i}_\varkappa$ after the $n$-th visit to $S_{i}$. Define the measure $\tilde{\nu}^{n,x,\varepsilon}_\varkappa$ on $\Gamma^{i}_\varkappa$ as follows:
\[
\tilde{\nu}^{n,x,\varepsilon}_\varkappa(B) = \mathrm{P}(A^{x,\varepsilon}_n \bigcap \{X^{x, \varepsilon}_{\tau^{x,\varepsilon}_n} \in B\} ),~~
B \in \mathcal{B}(\Gamma^{i}_\varkappa).
\]
Thus $\tilde{\nu}^{n,x,\varepsilon}_\varkappa$ is a sub-probability measure for each $n \geq 1$, while 
$\tilde{\nu}^{1 ,x,\varepsilon}_\varkappa = \tilde{\nu}^{x,\varepsilon}_\varkappa$ is a probability measure.
Observe that the process $X^{x,\varepsilon}_t$, starting at $x \in \Gamma^i_\varkappa$, reaches $S_i$ prior to reaching $S_j$ with probability
that is bounded from above by a constant $c < 1$, uniformly in $x$ and $\varepsilon$. (This is the case since $X^{x,\varepsilon}_t$ is
a perturbation of $X^x_t$, and $S_j$ is attracting.)
Therefore, 
\begin{equation} \label{exxx}
\tilde{\nu}^{n,x,\varepsilon}_\varkappa(\Gamma^i_\varkappa) =\mathrm{P}(A^{x,\varepsilon}_n) \leq
c^{n-1}~~{\rm for}~ n \geq 1.
\end{equation}

 Since $\lim_{\varepsilon \downarrow 0} f^\varepsilon(x)
= f(x)$, uniformly in $x \in \Gamma^{i}_\varkappa$, from Lemma~\ref{indmo} it follows that there is the limit
\[
\lim_{\varepsilon \downarrow 0} \int_{\Gamma^{i}_\varkappa} f^\varepsilon  d \tilde{\nu}^{1,x,\varepsilon}_\varkappa = 
\int_{\Gamma^{i}_\varkappa} f   d \tilde{\nu}_\varkappa =: r < 1. 
\]
By definition of $f^\varepsilon $ and $\tilde{\nu}^{1,x,\varepsilon}_\varkappa$, $r$ is the limiting  probability that the process
$X^{x,\varepsilon}_t$, starting at $x \in S_i$, returns to $S_i$ after one visit to $\Gamma^{i}_\varkappa$ prior to reaching $\partial D \setminus S_{i}$. By the strong Markov property (considering the first return to $S_i$), 
\[
\lim_{\varepsilon \downarrow 0} \int_{\Gamma^{i}_\varkappa} f^\varepsilon  d \tilde{\nu}^{2,x,\varepsilon}_\varkappa = 
 r^2 
\]
and, by induction,
\[
\lim_{\varepsilon \downarrow 0} \int_{\Gamma^{i}_\varkappa} f^\varepsilon  d \tilde{\nu}^{n,x,\varepsilon}_\varkappa = 
 r^n, ~~n \geq 1. 
\]
Therefore, by (\ref{exxx}), 
the measures $\sum_{n=1}^\infty \tilde{\nu}^{n,x,\varepsilon}_\varkappa$ converge weakly, uniformly in $x$,
to $(1+r+r^2+...) \tilde{\nu}_\varkappa$, that is, for each $h \in C(\Gamma^{i}_\varkappa)$, \be
\begin{equation} \label{wko}
\lim_{\varepsilon \downarrow 0} \int_{\Gamma^{i}_\varkappa}h d  (\sum_{n=1}^\infty \tilde{\nu}^{n,x,\varepsilon}_\varkappa) = 
\frac{1}{1-r}\int_{\Gamma^{i}_\varkappa} h d \tilde{\nu}_\varkappa,
\end{equation}
uniformly in $x \in S_{i}$. Since
\[
\mathrm{P}(X^{ x, \varepsilon}_{{\sigma}^{x, \varepsilon}} \in S_{j}) = \sum_{n = 1}^\infty
\int_{\Gamma^{i}_\varkappa} g^\varepsilon d  \tilde{\nu}^{n,x,\varepsilon}_\varkappa,
\]
from (\ref{wko}) and the uniform convergence of $g^\varepsilon$ to $g$, it follows that
\begin{equation} \label{rchh}
\lim_{\varepsilon \downarrow 0} Q^\varepsilon(x, S_{j}) = \lim_{\varepsilon \downarrow 0}  
\mathrm{P}(X^{ x, \varepsilon}_{{\sigma}^{x, \varepsilon}} \in S_{j}) = \frac{1}{1-r}\int_{\Gamma^{i}_\varkappa} g d \tilde{\nu}_\varkappa > 0,
\end{equation}
uniformly in $x \in S_{i}$. Although the integrand, the measure, and the set over which we integrate depend on $\varkappa$, the resulting integral in the right-hand side does not, since the left-hand side does not depend on $\varkappa$.

(b) Let $\eta > 0$. Fix $\varkappa > 0$ and $0 < \zeta < \varkappa$ such that (\ref{prex2aaz}) holds with $\rho  = \rho_i$, $\gamma = \gamma_i$, and
$S$, $\Gamma_\varkappa$, and $\Gamma_\zeta$ replaced by $S_{j}$, $\Gamma^{j}_\varkappa$, and $\Gamma^{j}_\zeta$, respectively. 
Since $S_j$ is repelling,  and $S_k$ is attracting for some $k \neq i,j$, by Lemma~\ref{lel1} (proved below) and Lemma~\ref{prat},     
we can make $\zeta$ smaller (if needed) so that
\begin{equation} \label{smpp}
\mathrm{P} (\tau^{x, \varepsilon}(\Gamma^j_\zeta) < \tau^{x,\varepsilon}(S_k)) \leq \eta
\end{equation}
for all $x \in \Gamma^j_\varkappa$ and all sufficiently small $\varepsilon$.
For $x \in S_i$, let
\[
\tilde{\sigma}^{{x},\varepsilon} = \inf\{t: {{X}}^{ x,\varepsilon}_t \in (\partial D \bigcup \Gamma^j_\zeta )\setminus S_i\}.
\]
We modify the definition of the functions $f$ and $g$ from part (a) of the proof. Now
\[
f(x) = \mathrm{P}\left(\lim_{t \rightarrow \infty} {\rm dist}(X^x_t, S_{i})  = 0~~{\rm and}~~X^x_t \notin \Gamma^j_\zeta~~{\rm for}~t \geq 0\right),
\]
\[
g(x) = \mathrm{P}\left(X^x_t \in \Gamma^j_\zeta~~{\rm for}~{\rm some}~t > 0 \right).
\]
The same arguments that led to (\ref{rchh}) now give
\begin{equation} \label{rchh2}
\lim_{\varepsilon \downarrow 0}  
\mathrm{P}(X^{ x, \varepsilon}_{\tilde{\sigma}^{x, \varepsilon}} \in \Gamma^j_\zeta) = \frac{1}{1-r}\int_{\Gamma^{i}_\varkappa} g d \tilde{\nu}_\varkappa > 0,
\end{equation}
uniformly in $x \in S_i$. Therefore, from the strong Markov property and the lower bound from Lemma~\ref{prre}, we get the lower bound:
\begin{equation} \label{lobd}
\mathrm{P}(X^{ x, \varepsilon}_{{\sigma}^{x, \varepsilon}} \in S_j) \geq 
\rho_j (1 - \eta)^2  \frac{ \zeta^{\gamma_j}}{\varepsilon^{\gamma_j}} \frac{1}{1-r}\int_{\Gamma^{i}_\varkappa} g d \tilde{\nu}_\varkappa 
\end{equation}
for all sufficiently small $\varepsilon$. The extra factor $(1-\eta)$ in the right-hand side is due to the fact that the equality 
in (\ref{rchh2}) is achieved only for the limit. In order to get an upper bound, we take into account the contributions to 
$\mathrm{P}(X^{ x, \varepsilon}_{{\sigma}^{x, \varepsilon}} \in S_j)$ from the following events:

$E_0$:  $X^{x,\varepsilon}_t$ reaches $\Gamma^j_\zeta$ and proceeds to $S_j$ before reaching $\Gamma^j_\varkappa$,

$E_1$: $X^{x,\varepsilon}_t$ reaches $\Gamma^j_\zeta$, and proceeds to $S_j$ after one excursion to $\Gamma^j_\varkappa$ and back
to  $\Gamma^j_\zeta$, but  before reaching $S_k$.

$E_n$: the same, but with $n$ excursions, $n \geq 2$. 
\\
From the strong Markov property, (\ref{smpp}), (\ref{rchh2}), and the upper bound from Lemma~\ref{prre}, we get 
\[
\mathrm{P}(E_n) \leq 
\rho_j (1 + \eta)^2 \eta^n  \frac{ \zeta^{\gamma_j}}{\varepsilon^{\gamma_j}} \frac{1}{1-r}\int_{\Gamma^{i}_\varkappa} g 
d \tilde{\nu}_\varkappa 
\]
for all $x \in S_i$, all $n$, and all sufficiently small $\varepsilon$.
Taking the sum in $n$ gives us the upper bound:
\begin{equation} \label{upbd}
\mathrm{P}(X^{ x, \varepsilon}_{{\sigma}^{x, \varepsilon}} \in S_j) \leq \rho_j \frac{(1 + \eta)^2}{1-\eta}    \frac{ \zeta^{\gamma_j}}{\varepsilon^{\gamma_j}} \frac{1}{1-r}\int_{\Gamma^{i}_\varkappa} g 
d \tilde{\nu}_\varkappa 
\end{equation}
for all $x \in S_i$  and all sufficiently small $\varepsilon$. Combining (\ref{lobd}) and (\ref{upbd}), we obtain
\[
(1-\eta)^2 G(\varkappa, \zeta) \leq \varepsilon^{\gamma_{j}} Q^\varepsilon(x, S_{j}) \leq \frac{(1 + \eta)^2}{1-\eta}
G(\varkappa, \zeta),
\]
where $G(\varkappa, \zeta) = \rho_j    { \zeta^{\gamma_j}}  ({1-r})^{-1}\int_{\Gamma^{i}_\varkappa} g 
d \tilde{\nu}_\varkappa$. Since $\eta$ can be taken arbitrarily small (by selecting appropriate $\varkappa$ and $\zeta$), and $Q^\varepsilon(x, S_{j})$ does not depend on $\varkappa$ and $\zeta$, there is a positive limit
\[
\lim_{\varepsilon \downarrow 0} \varepsilon^{\gamma_{j}} Q^\varepsilon(x, S_{j}) = : q_{i j},
\]
uniformly in $x \in S_i$. 

(c) We provide only a sketch of a proof here in order to avoid unnecessary technical details. If there are fewer than two attracting components of
the boundary, we take two disjoint $(d-1)$-dimensional spheres of radius $\delta > 0$ centered at $x_1, x_2 \in D$ 
and denote them by $S_{m+1}$ and $S_{m+2}$.  
The definitions of the sequence of stopping times ${\sigma}^{{x},\varepsilon}_n$ and of the corresponding kernel $Q^\varepsilon$ need  to 
be modified to allow transitions to and from $S_{m+1}$ and $S_{m+2}$. These two surfaces are designated as attracting, although they are
not invariant for the unperturbed process. We have the following modifications of parts (a) and (b) of the lemma: for each $\eta > 0$,
sufficiently small radius $\delta$ can be chosen for the spheres so that $
|Q^\varepsilon(x, S_{j})  - q_{i j}| < \eta$, $x \in S_{i}$ if $S_j$ is attracting and
$|\varepsilon^{\gamma_{j}} Q^\varepsilon(x, S_{j}) - q_{i j}| < \eta$, $x \in S_{i}$ if $S_j$ is repelling. We are interested not
in $Q^\varepsilon(x, S_j)$ now, but in the asymptotics of the event that $X^{x,\varepsilon}_t$ starting at $x \in S_i$ ($ i \leq m$) reaches
$S_j$ prior to reaching $\partial D \setminus (S_i \bigcup S_j)$ (i.e., visits to $S_{m+1}$ and $S_{m+2}$ are allowed). The asymptotics of 
all the transition probabilities are described, though (up to a small extra term $\eta$), and therefore the desired result easily follows. 
\qed
\\
\\
{REMARK.} In the case (a), when $S_j$ is attracting, the limiting probability can be described as follows. For $x \in D$, 
let $\mathbf{X}^{i, x}_t$ be the process that is obtained from $X^x_t$ by conditioning on the complement of the event that $\lim_{t \rightarrow \infty}  {\rm dist}(X^x_t, S_k) = 0$. 
Let $\mathbf{p}^{i,x}_j = \mathrm{P}(\lim_{t \rightarrow \infty} {\rm dist} (\mathbf{X}^{i, x}_t, S_j) = 0)$. It is possible to
show that $q_{ij} = \lim_{{\rm dist}(x, S_i) \downarrow 0} \mathbf{p}^{i,x}_j$. 

\section{Metastable distributions for the perturbed process} \label{mdis}

In this section, we describe the distribution of the process $X^{x,\varepsilon}_t$ and the stopped process $X^{x,\varepsilon}_{t \wedge 
\tau^{x, \varepsilon}(\partial D)}$ at different time scales. Recall that, by Lemma~\ref{spec}, there is a constant $\gamma_k$ associated to each
component of the boundary $S_k$. Without loss of generality, we assume that $\gamma_1 \geq \gamma_2 \geq ... \geq \gamma_m$. 

\subsection{The case when at least one  component of the boundary is attracting}

Here, we assume that $\gamma_1 > 0$, i.e., $S_1$ is attracting. Let $\overline{m}$ be such that
$\gamma_1 \geq ... \geq \gamma_{\overline{m}} > 0 > \gamma_{\overline{m}+1} \geq ... \geq \gamma_m$. For $x \in D$ and $k \leq \overline{m}$, let $E^x_k$ be the event that $\lim_{t \rightarrow \infty}  {\rm dist}(X^x_t, S_k) = 0$.  
\begin{lemma} \label{lel1}
For each $x \in D$, $\mathrm{P}(\bigcup_{k =1}^{\overline{m}} E^x_k )= 1$.
\end{lemma}
\proof
Let $\eta > 0$. Recall the definition of the sets $V^k_{\varkappa_1, \varkappa_2}$ from the beginning of Section~\ref{neigh} (we
now use the superscript $k$ to refer to a particular component of the boundary). Let
$\overline{U}_\varkappa$ be the closure of the set $ D \setminus(\bigcup_{k = 1}^m V^k_{0,\varkappa})$. Thus 
$\overline{U}_\varkappa$ is a closed set obtained by
removing a small neighborhood of the boundary from $D$.  

Take $\varkappa > 0$ and $0 < \zeta < \varkappa$ sufficiently small so that (\ref{prex2aa}) holds (for all sufficiently small $\varepsilon$) for each attracting component of the boundary with the right-hand side that admits the estimate $ (1 + \eta)
 { \zeta^\gamma}/{\varkappa^\gamma} < \eta$. Since $X^{x,\varepsilon}_t$ converges to $X^x_t$ on each finite time interval, 
we conclude that, conditioned on reaching $\Gamma^k_\zeta$ with $k \leq \overline{m}$, the process $X^x_t$ forever remains in $V^k_{0,\varkappa}$ with probability at least $1 - \eta$. Since arbitrarily small $\eta$ and $\varkappa$ can be taken, the result will follow if
we show that  $\mathrm{P}(\tau^x(\bigcup_{k=1}^{\overline{m}} \Gamma^k_\zeta) < \infty) = 1$, where
$\tau^{x}(A) = \inf\{t \geq 0: X^{x}_t \in A\}$ for a closed set $A$.

Since $X^x_t$ is non-degenerate away from the boundary, there is $p = p(\varkappa, \zeta)> 0$ such that
$\mathrm{P}(\tau^x(\bigcup_{k=1}^{\overline{m}} \Gamma^k_\zeta) < \tau^x(\bigcup_{k=\overline{m}+1}^{m} \Gamma^k_\zeta)) \geq p$ 
 for all $x \in \overline{U}_\varkappa$.
On the other hand, it follows from  Lemma~\ref{prre} that $\mathrm{P}(\tau^x(\overline{U}_\varkappa) < \infty) = 1$ for $x \in \bigcup_{k=\overline{m}+1}^{m} \Gamma^k_\zeta$. Therefore, from the strong Markov property, it follows that $X^x_t$ will reach
$\bigcup_{k=1}^{\overline{m}} \Gamma^k_\zeta$ after a finite number of excursions between $\overline{U}_\varkappa$ and $\bigcup_{k=\overline{m}+1}^{m} \Gamma^k_\zeta$. 
\qed
\\

Now we can describe the distribution of the perturbed process  at time scales $t(\varepsilon)$ that satisfy $1 \ll t(\varepsilon) \ll 
\varepsilon^{-\gamma_{\overline{m}}}$.  Recall that $\pi_k$ are the invariant measures for the process $X^x_t$ considered as a process on $S_k$. Let $p^x_k = \mathrm{P}(E^x_k)$, $1 \leq k \leq \overline{m}$. 

\begin{theorem} \label{iuo} If $1 \ll t(\varepsilon) \ll \varepsilon^{-\gamma_{\overline{m}}}$ and $x \in D$, then the distribution of 
$X^{x,\varepsilon}_{t(\varepsilon)}$ converges to the measure $\sum_{k=1}^{\overline{m}} p^x_k \pi_k$. 
\end{theorem}
\proof 
 From Lemmas~\ref{lel1}, \ref{prat}, and \ref{mltime} (part (c)), and the proximity of $X^{x,\varepsilon}_t$ and $X^x_t$ on finite
time intervals, it follows that, for each $\zeta > 0$,
\[
\lim_{\varepsilon \downarrow 0} \mathrm{P}(X^{x,\varepsilon}_{t(\varepsilon)/2} \in V^k_{0,\zeta}) = p^x_k
\]
for each $1 \leq k \leq \overline{m}$. The result now follows from Lemma~\ref{wcon}.
\qed
\\

Next, let us explore the behavior of $X^{x,\varepsilon}_t$ at longer time scales. Assume that $\varepsilon^{-\gamma_{l+1}} \ll t(\varepsilon)
\ll \varepsilon^{-\gamma_{l}}$, where $l+1 \leq \overline{m}$ and $l \geq 1$ (this is possible if $\overline{m} \geq 2$).  
Consider the  discrete-time 
Markov chain $Z_n^x$ on $\{1,...,\overline{m}\}$ with transition probabilities $q_{ij}$ defined in Part (a) of Lemma~\ref{trapl} (it follows from
Lemma~\ref{trapl} that these do form a stochastic matrix). We take the vector $(p^x_1,...,p^x_{\overline{m}})$ as the initial distribution
of the Markov chain. Let $\tau_l$ be the hitting time of the set $\{1,...l\}$, and define
\[
p^{x,l}_k = \mathrm{P}(Z^x_{\tau_l} = k),~~~k =1,...,l.
\]
\begin{theorem} \label{btto} If $\varepsilon^{-\gamma_{l+1}} \ll t(\varepsilon)
\ll \varepsilon^{-\gamma_{l}}$ and $x \in D$, then the distribution of 
$X^{x,\varepsilon}_{t(\varepsilon)}$ converges to the measure $\sum_{k=1}^{l} p^{x,l}_k \pi_k$. 
\end{theorem}
\proof From Lemmas~\ref{lel1} and \ref{prat}, it follows that $\mathrm{P}(X^{x,\varepsilon}_{\tau^{x,\varepsilon} (\partial D)} \in S_k) \rightarrow 
p^x_k$, as $\varepsilon \downarrow 0$, for $1 \leq k \leq {\overline{m}}$. From Lemma~\ref{trapl} it then follows that 
$\mathrm{P}(X^{x,\varepsilon}_{\tau^{x,\varepsilon}(S_1 \bigcup ... \bigcup S_l)} \in S_k) \rightarrow 
p^{x,l}_k$, as $\varepsilon \downarrow 0$, for $1 \leq k \leq l$. The set $S_1 \bigcup ... \bigcup S_l$ is reached after the first transition
from $x$ to $\partial D$ (which takes  time $O(|\ln(\varepsilon)|)$ due to Lemmas~\ref{lel1}, \ref{lee2} (part (a)), and \ref{trp5}
 (part (a))) and
a finite number of transitions between different components of the boundary, which take time $O(\varepsilon^{-\gamma_{l+1}})$ by 
Lemma~\ref{mltime} (part (a)). Therefore,
\[
\lim_{\varepsilon \downarrow 0} \mathrm{P} ( \tau^{x,\varepsilon}(S_1 \bigcup ... \bigcup S_l) < \frac{t(\varepsilon)}{2})  = 1.
\]
By Lemma~\ref{mltime} (part (c)) and the strong Markov property, for each $\zeta > 0$, 
\[
\lim_{\varepsilon \downarrow 0} \mathrm{P} ( X^{x,\varepsilon}_{t(\varepsilon)/2} \in V^k_{0,\zeta}) = p^{x,l}_k
\]
for each $1 \leq k \leq l$. The result now follows from Lemma~\ref{wcon}.
\qed
\\

Next, we explore the longest time scales. Here, for simplicity, we assume that $\gamma_1 > \gamma_2$. 
\begin{theorem} \label{ltsc} If $\gamma_1 > \gamma_2$, 
$ t(\varepsilon) \gg \varepsilon^{-\gamma_{1}}$,
 and $x \in D$, then the distribution of 
$X^{x,\varepsilon}_{t(\varepsilon)}$ converges to the measure $\pi_1$. 
\end{theorem}
\proof Take an arbitrary function $t_1(\varepsilon)$ that satisfies $\varepsilon^{-\gamma'} \ll t_1(\varepsilon)
\ll \varepsilon^{-\gamma_{1}}$ with some $\gamma'$ such that $\max(0, \gamma_2) < \gamma' < \gamma_1$. Let $\mu^{x,\varepsilon}$ be
the distribution of $X^{x,\varepsilon}_{t(\varepsilon) - t_1(\varepsilon)}$.  Then the distribution of $X^{x,\varepsilon}_{t(\varepsilon)}$
agrees with that of  $X^{\mu^{x,\varepsilon},\varepsilon}_{t_1(\varepsilon)}$ (the process whose initial distribution is $\mu^{x,\varepsilon}$
rather than concentrated in one point). The same proof as for Theorem~\ref{iuo} (if $\gamma_2 < 0$) or Theorem~\ref{btto} (if
$\gamma_2 > 0$) now applies; the only difference is that we now need to bound the time $\tau^{\mu^{x,\varepsilon} ,\varepsilon}(\partial D)$.
From Lemmas~\ref{lee1} (part (a)), \ref{lee2} (part (a)), and \ref{trp5} (part (a)), it easily follows that, for each $\eta > 0$, there is $C > 0$ such that
\[
\mathrm{P}(\tau^{x, \varepsilon}(\partial D) < C |\ln(\varepsilon)|) \geq 1 -\eta
\]
for all sufficiently small $\varepsilon$ uniformly in   $x \in \overline{D}$. This implies the desired result.
\qed
\\
\\
{REMARK.} At the `transitional' time scales ($t(\varepsilon) \sim \varepsilon^{-\gamma_l}$ with some $l \leq \overline{m}$) the
limiting distribution will be a linear combination of the measures $\pi_k$, $k \leq l$. Specifying the coefficients in this linear combination
requires a yet more delicate analysis of the time it takes the process $X^{x,\varepsilon}_t$ to exit  a neighborhood of an attracting surface, and is not done here. 
\\

We can also describe the distribution of the stopped process. Let $\nu_k$ be the measure $\nu$ defined in Lemma~\ref{mlee1} (with $S_k$ instead of a generic component $S$).
\begin{theorem}  
\label{iuox} If $1 \ll t(\varepsilon) \ll |\ln(\varepsilon)|$ and $x \in D$, then the distribution of 
 $X^{x,\varepsilon}_{t(\varepsilon) \wedge 
\tau^{x, \varepsilon}(\partial D)}$ converges to the measure $\sum_{k=1}^{\overline{m}} p^x_k \pi_k$. If 
$ t(\varepsilon) \gg |\ln(\varepsilon)|$ and $x \in D$, then the distribution of 
 $X^{x,\varepsilon}_{t(\varepsilon) \wedge 
\tau^{x, \varepsilon}(\partial D)}$ converges to the measure $\sum_{k=1}^{\overline{m}} p^x_k \nu_k$.
\end{theorem}
\proof If $1 \ll t(\varepsilon) \ll |\ln(\varepsilon)|$, then 
$\lim_{\varepsilon \downarrow 0} \mathrm{P} (\tau^{x, \varepsilon}(\partial D) < t(\varepsilon)) = 0$,
as follows from Lemmas~\ref{lee2} (part (b)) and \ref{lel1}. Thus the result follows from Theorem~\ref{iuo}. 
If $ t(\varepsilon) \gg |\ln(\varepsilon)|$, then 
$\lim_{\varepsilon \downarrow 0} \mathrm{P} (\tau^{x, \varepsilon}(\partial D) < t(\varepsilon)) = 1$, and the result
follows from Lemma~\ref{mlee1}.
\qed

\subsection{The case when all the  components of the boundary are repelling}
Here, we assume that $\gamma_1 < 0$, i.e., all the components of the boundary are repelling.  
\begin{theorem} \label{hujuh} There is a unique measure $\mu$ on $D$ that is invariant for the process $X^x_t$, $x \in D$. This measure
is absolutely continuous with respect to the Lebesgue measure. If $t(\varepsilon) \rightarrow \infty$ as $\varepsilon \downarrow 0$, then
the distribution of $X^{x,\varepsilon}_{t(\varepsilon)}$ converges to $\mu$ for each $x \in D$. 
\end{theorem}
\proof Consider two $(d-1)$-dimensional spheres $F$ and $G$ such that $F \subset D$ and $G$ is contained in 
the interior of the ball bounded by $F$. Let 
\[
\sigma^{x,\varepsilon}_F = \inf \{t \geq 0: X^{x,\varepsilon}_t \in F \}, ~~~
\sigma^{x }_F = \inf \{t \geq 0: X^{x }_t \in F \},
\]
\[
\sigma^{x,\varepsilon}_G = \inf \{t \geq 0: X^{x,\varepsilon}_t \in G \},~~~
\sigma^{x }_G = \inf \{t \geq 0: X^{x }_t \in G \}.
\]

Consider  $\zeta$ and $\varkappa$ such that $0 < \zeta < \varkappa$.
By considering successive visits  of the process $X^{x,\varepsilon}_t$, prior to reaching $F$,  to the sets $\bigcup_{k =1}^n \Gamma^k_\zeta$ 
and  $\bigcup_{k =1}^n \Gamma^k_\varkappa \bigcup S$, and using Lemmas~\ref{lee1} (parts (a) and (c)) and \ref{prre},  we see that 
the random variables $\sigma^{x,\varepsilon}_G $ are uniformly integrable in $\varepsilon$ and $x \in F$.

From the non-degeneracy of the process in $D$, it follows that $ \sigma^{x,\varepsilon}_F $ are uniformly integrable in $\varepsilon$ 
and $x \in G$. 

For $x \in G$, let
\[
\tilde{\sigma}^{x,\varepsilon} = \inf\{t \geq \sigma^{x, \varepsilon}_F: X^{x,\varepsilon}_t \in G\},~~~
\tilde{\sigma}^{x } = \inf\{t \geq \sigma^{x }_F: X^{x }_t \in G\}.
\]
Thus  $\tilde{\sigma}^{x,\varepsilon}$ ($\tilde{\sigma}^{x }$) is the first time when the corresponding process returns 
to $G$ after visiting $F$. From the proximity of $X^{x,\varepsilon}_t$ and $X^{x}_t$, we get that $\tilde{\sigma}^{x,\varepsilon} \rightarrow \tilde{\sigma}^{x }$ almost surely as $\varepsilon \downarrow 0$.  The uniform integrability of 
$\tilde{\sigma}^{x,\varepsilon}$ (which follows from the uniform integrability of  $\sigma^{x,\varepsilon}_G $, $x \in F$, $\varepsilon > 
0$,  and 
$\sigma^{x,\varepsilon}_F $, $x \in G$, $\varepsilon > 
0$)  implies 
that $\tilde{\sigma}^{x}$ are uniformly integrable. 

We can introduce Markov kernels $Q^\varepsilon $ and $Q $ on $G$ via
\[
Q^\varepsilon(x, A)  = \mathrm{P}(X^{x,\varepsilon}_{\tilde{\sigma}^{x,\varepsilon}} \in A),~~Q (x, A)  = 
\mathrm{P}(X^{x }_{\tilde{\sigma}^{x }} \in A),~~~x \in G,~~A \in \mathcal{B}(G).
\]
Since the processes $X^{x,\varepsilon}_t$ and $X^{x}_t$ are non-degenerate in a neighborhood of $G$, $Q^\varepsilon(x,\cdot)$ and
$Q(x, \cdot)$ have densities that are uniformly bounded from above and below. Therefore, there exist unique probability measures $P^\varepsilon$ and $P$ on $G$ such that
\[
P^\varepsilon Q^\varepsilon = P^\varepsilon,~~PQ = P.
\]
Since the expectation of $\tilde{\sigma}^{x,\varepsilon}$ ($\tilde{\sigma}^{x }$) is bounded, the invariant measure $\mu^\varepsilon$ ($\mu$) for the process $X^{x,\varepsilon}_t$ ($X^{x}_t$) on $\overline{D}$ ($D$) can
now be expressed explicitly
\begin{equation} \label{exexp}
\mu^\varepsilon (A) = 
\frac{\int_G \mathrm{E} \int_0^{\tilde{\sigma}^{x,\varepsilon}} \chi_A (X^{x,\varepsilon}_t) dt d P^\varepsilon(x)}{ \int_G \mathrm{E}
\tilde{\sigma}^{x,\varepsilon} d P^\varepsilon(x)},~~~~
\mu (A) = 
\frac{\int_G \mathrm{E} \int_0^{\tilde{\sigma}^{x }} \chi_A (X^{x }_t) dt d P (x)}{ \int_G \mathrm{E}
\tilde{\sigma}^{x } d P (x)},~~~A \in \mathcal{B}(D).
\end{equation}
(Note that  $\mu^\varepsilon(\partial D) =0$, and thus  $\mu^\varepsilon$ is a probability measure on $D$. ) From the
proximity of $X^{x,\varepsilon}_t$ and $X^x_t$ on finite time intervals it follows that $P^\varepsilon$ converges weakly
to $P$ as $\varepsilon \downarrow 0$. 

By (\ref{exexp}), the uniform integrability of $\tilde{\sigma}^{x,\varepsilon}$,  and due to the
proximity of $X^{x,\varepsilon}_t$ and $X^x_t$ on finite time intervals, $\mu^\varepsilon$ converges to $\mu$. It remains to note that, for
each $f \in C_b(D)$, 
\[
\lim_{t \rightarrow \infty}(\mathrm{E} f(X^{x,\varepsilon}_t) - \int_D f d \mu^\varepsilon) = 0  
\]
uniformly in $\varepsilon > 0$. Therefore, 
the distribution of $X^{x,\varepsilon}_{t(\varepsilon)}$ converges to $\mu$ if $t(\varepsilon) \rightarrow \infty$. 
\qed
\\

We can also describe the distribution of the stopped process. Here we assume, for simplicity, that $\lambda_1 > \lambda_2$, i.e., $S_1$ 
is the ``least repelling" component of the boundary.
\begin{theorem}
\label{iuo22} If $1 \ll t(\varepsilon)     \ll \varepsilon^{-\gamma_1}$ and $x \in D$, then the distribution of 
 $X^{x,\varepsilon}_{t(\varepsilon) \wedge 
\tau^{x, \varepsilon}(\partial D)}$ converges to the measure $\mu$. If 
$ t(\varepsilon) \gg  \varepsilon^{-\gamma_1}$ and $x \in D$, then the distribution of 
 $X^{x,\varepsilon}_{t(\varepsilon) \wedge 
\tau^{x, \varepsilon}(\partial D)}$ converges to the measure $ \nu_1$.
\end{theorem}  
\proof Using  Lemma \ref{prre} together with the strong Markov property of the process, it is not difficult to show that
\[
\lim_{\varepsilon \downarrow 0} \mathrm{P}(\tau^{x, \varepsilon}(\partial D) < t(\varepsilon)) = 0,~~~{\rm if}~~~t(\varepsilon) \ll \varepsilon^{-\gamma_1},
\]
\[
\lim_{\varepsilon \downarrow 0} \mathrm{P}\left(\tau^{x, \varepsilon}(\partial D) < t(\varepsilon),~~X^{x,\varepsilon}_{\tau^{x, \varepsilon}(\partial D)} \in S_1\right) = 1,~~~
{\rm if}~~~t(\varepsilon) \gg \varepsilon^{-\gamma_1}.
\]
Thus the first statement of the theorem follows from Theorem~\ref{hujuh}, while the second part follows from Lemma~\ref{mlee1}
\qed
\\

By the stochastic representation of solutions to parabolic PDEs, Theorem~\ref{mnt1} follows from Theorems~\ref{iuox}
and \ref{iuo22}. Similarly, Theorem~\ref{mnt2} follows from Theorems~\ref{iuo}, \ref{btto},  \ref{ltsc}, and \ref{hujuh}.

\section{Periodic homogenization} \label{homog1}
In this section, we consider  processes (and the corresponding operators) in the entire space $\mathbb{R}^d$ rather than in 
a bounded domain. Assume that the coefficients of the operators $L$ and $\tilde{L}$ defined in Section~\ref{intro} are one-periodic in each 
of the variables. We will assume that the coefficients of $L$ degenerate (as in Section~\ref{intro},  in the direction normal to the
surface) on a periodic array of $C^4$-surfaces $S_{l,k}$, $ l \in \mathbb{Z}^d$, $1 \leq k \leq m$, that serve as boundaries of bounded domains
$D_{l,k}$. The domains themselves are assumed to be disjoint: $\overline{D}_{l_1, k_1} \bigcap \overline{D}_{l_2,k_2} = \emptyset$ if
$(l_1, k_1) \neq (l_2, k_2)$, and therefore the complement $D = \mathbb{R}^d \setminus (\bigcup_{l, k} D_{l,k})$ is an unbounded connected set.

As before, each of the surfaces $S_{l,k}$ can be classified as attracting or repelling depending on the sign of $\gamma_{l,k}$ 
(see Lemma~\ref{spec}). By periodicity, there are constants $\gamma_k$, $1 \leq k \leq m$, such that $\gamma_k = \gamma_{l,k}$ for each $l \in \mathbb{Z}^d$.
As in Section~\ref{mdis}, we assume, without loss of generality, that $\gamma_1 \geq \gamma_2 \geq ... \geq \gamma_m$. 

In this section, $X^{x,\varepsilon}_t$  again stands for the process satisfying (\ref{peqn1}), however, now we assume that there
is no reflection on the surfaces $S_{l,k}$, and thus the process takes values in the entire space $\mathbb{R}^d$. 
One can expect that, due to the invariance of the coefficients with respect to integer shifts, the process $X^{x,\varepsilon}_t$ can be
approximated, in appropriate space-time scales, by a diffusion process with constant drift and diffusion coefficients. Results of this type are considered in this section. We will consider two cases: when either all the surfaces are attracting or when all the surfaces are
repelling; the general case can be studies using a combination of these two scenarios. We only provide sketches of the proofs here
since the arguments are largely similar to those used in the earlier sections. 


We start with the case when $\gamma_m > 0$, i.e., all the surfaces are attracting. In this case, at time scales larger 
than $\varepsilon^{-\gamma_1}$, the
behavior of the process $X^{x,\varepsilon}_t$ is diffusive. Namely, we have the following result.
\begin{theorem} If $\gamma_m >0$ and $t(\varepsilon) \gg \varepsilon^{-\gamma_1}$, then there is  a vector $a \in \mathbb{R}^d$ such that, for any $x \in \mathbb{R}^d$,
\[
\lim_{\varepsilon \downarrow 0}  \frac{ \mathrm{E} X^{x,\varepsilon}_{t(\varepsilon)}}{\varepsilon^{\gamma_1} 
t(\varepsilon)}  = a;~~~ \lim_{\varepsilon \downarrow 0} \frac{   X^{x,\varepsilon}_{t(\varepsilon)}}{\varepsilon^{\gamma_1} 
t(\varepsilon)}  = a~~~in~probability.
\]
Moreover, there is a $d\times d$ positive-definite symmetric matrix $B$
such that
\[
\frac{X^{x,\varepsilon}_{t(\varepsilon)} -\mathrm{E} X^{x,\varepsilon}_{t(\varepsilon)}}{\sqrt{ \varepsilon^{\gamma_1} t(\varepsilon)}} \rightarrow N(0, B)~~~in~~distribution,~~
as~~\varepsilon \downarrow 0. 
\]
\end{theorem}
\noindent
{\it Sketch of the proof.} Assume that $x \in \partial D$ (the general case is treated similarly, since $X^{x,\varepsilon}_t$ reaches
$\partial D$ sufficiently fast). As in Section~\ref{atpr}, we can consider a Markov renewal process
$({X}^{ x, \varepsilon}_{{\sigma}^{{x, \varepsilon}}_{n}},  
{\sigma}^{{x, \varepsilon}}_{n})$, $n \geq 0$,
 on $\partial D$ (which is now a periodic
array of surfaces). In particular, $Y^{x, \varepsilon}_n := {X}^{ x, \varepsilon}_{{\sigma}^{{x, \varepsilon}}_{n}}$ is 
a discrete-time Markov chain on $\partial D$.

By the Law of Large Numbers  for Markov chains, there exist a vector $\tilde{a}(\varepsilon)$
such that
\begin{equation} \label{lln1}
\lim_{n \rightarrow 0}  \frac{ \mathrm{E} Y^{x,\varepsilon}_{n}}{n}  = \tilde{a}(\varepsilon);~~~ \lim_{n \rightarrow 0}  
\frac{   Y^{x,\varepsilon}_{n}}{n}  = \tilde{a}(\varepsilon)~~~{\rm in}~{\rm probability}.
\end{equation}
By the Central Limit Theorem for Markov chains,
there exists
a matrix $\tilde{B}(\varepsilon)$ such that
\begin{equation} \label{cltd}
\frac{Y^{x,\varepsilon}_n -  \mathrm{E} Y^{x,\varepsilon}_n }{\sqrt{ n}} 
\rightarrow N(0, \tilde{B}(\varepsilon))~~~{\rm in}~~{\rm distribution},~~
{\rm as}~~n \rightarrow \infty. 
\end{equation}
The behavior, as $\varepsilon \downarrow 0$, of the transition kernel $Q^\varepsilon$ of $Y^{x,\varepsilon}_n$ 
can be understood using Lemma~\ref{trapl},
which is still applicable (despite the process 
$X^{x,\varepsilon}_t$ being considered in the entire space and not reflected on the boundary of $D$). 
Namely, there   exist constants $q_{(l_1, k_1), (l_2,k_2)} > 0$, $l_1, l_2 \in \mathbb{Z}^d$, $1 \leq k_1, k_2 \leq m$, 
$(l_1, k_1) \neq (l_2,k_2)$, such that 
\[
\lim_{\varepsilon \downarrow 0} Q^\varepsilon(x, S_{l_2,k_2}) = q_{(l_1, k_1), (l_2,k_2)},~~x \in S_{l_1, k_1}.
\]
Thus, for small $\varepsilon$, the long-time behavior of $Y^{x,\varepsilon}_n$ is close to that of a spatially homogeneous random 
walk on the lattice $\mathbb{Z}^d \times \{1,...,m\}$ with transition probabilities $q_{(l_1, k_1), (l_2,k_2)}$. The latter do  not
depend on  $\varepsilon$, and the Law of Large Numbers and Central Limit Theorem with a drift $\tilde{a}$ and a diffusion matrix $\tilde{B}$ applies to the random walk. 
Thus $\lim_{\varepsilon \downarrow 0} \tilde{a}(\varepsilon) = \tilde{a}$ and 
$\lim_{\varepsilon \downarrow 0} \tilde{B}(\varepsilon) = \tilde{B}$.

It is not difficult to show that the limiting transition kernel  $Q^\varepsilon$ has the following properties, uniformly in $\varepsilon$:

(a) For each $C > 0$, there is $c > 0$ such that $Q^\varepsilon(x, S_{l_2,k_2}) \geq c$, provided that $x \in S_{l_1, k_1}$ with $|l_1 - l_2| 
\leq C$ (indeed, due to non-degeneracy of the process in $D$, the probability of a transition between nearby surfaces is bounded
from below);

(b) There is $\lambda > 0$ such that $Q^\varepsilon(x, S_{l_2,k_2}) \leq e^{-\lambda |l_1 - l_2|}$, provided that $x \in S_{l_1, k_1}$
(the probability that $X^{x,\varepsilon}_t$ travels far without hitting
any of the attracting surfaces is small).

From these two properties, it follows that the convergence in (\ref{lln1}) and (\ref{cltd}) is uniform with respect to 
$\varepsilon$ in the sense that, for each 
integer-valued function $n(\varepsilon)$ such that $\lim_{\varepsilon \downarrow 0} 
n(\varepsilon)= \infty$, 
\begin{equation} \label{lln1a}
\lim_{\varepsilon \downarrow 0}  \frac{ \mathrm{E} Y^{x,\varepsilon}_{n(\varepsilon)}}{n(\varepsilon)}  = 
\tilde{a};~~~ \lim_{\varepsilon \downarrow 0}   
\frac{   Y^{x,\varepsilon}_{n(\varepsilon)}}{n(\varepsilon)}  = \tilde{a}~~{\rm in}~{\rm probability},
\end{equation}
\begin{equation} \label{cltda}
\frac{Y^{x,\varepsilon}_{n(\varepsilon)} -  \mathrm{E} Y^{x,\varepsilon}_{n(\varepsilon)} }{\sqrt{ n(\varepsilon)}} 
\rightarrow N(0, \tilde{B} )~~~{\rm in}~~{\rm distribution},~~
{\rm as}~~ {\varepsilon \downarrow 0} . 
\end{equation}
Moreover, (\ref{lln1a}) and (\ref{cltda}) can be modified to allow the time at which the process is evaluated to be a small random perturbation of a deterministic quantity. Namely,
\begin{equation} \label{lln1b}
\lim_{\varepsilon \downarrow 0}  \frac{ \mathrm{E} Y^{x,\varepsilon}_{N(\varepsilon)}}{n(\varepsilon)}  = 
\tilde{a};~~~ \lim_{\varepsilon \downarrow 0}   
\frac{   Y^{x,\varepsilon}_{N(\varepsilon)}}{n(\varepsilon)}  = \tilde{a}~~{\rm in}~{\rm probability},
\end{equation}
\begin{equation} \label{cltdb}
\frac{Y^{x,\varepsilon}_{N(\varepsilon)} -  \mathrm{E} Y^{x,\varepsilon}_{N(\varepsilon)} }{\sqrt{ n(\varepsilon)}} 
\rightarrow N(0, \tilde{B} )~~~{\rm in}~~{\rm distribution},~~
{\rm as}~~ {\varepsilon \downarrow 0}, 
\end{equation}
where $N(\varepsilon) = n(\varepsilon) + \xi(\varepsilon)$ with $\xi(\varepsilon)/ {n(\varepsilon)} \rightarrow 0$ in $L^1$ as $\varepsilon \downarrow 0$.  
The validity of  (\ref{lln1b}) and (\ref{cltdb})  can be proved using slightly 
stronger versions of (\ref{lln1}) and (\ref{cltd}) (e.g., the invariance principle instead of the CLT for Markov chains). 
 
To complete the proof of the theorem, we need to apply (\ref{lln1a}) and (\ref{cltda}) with $N(\varepsilon) = 
\max\{n: {\sigma}^{{x, \varepsilon}}_{n} \leq t(\varepsilon) \}$. We claim that there are a constant $c > 0$  such that
\begin{equation} \label{annk}
N(\varepsilon) = c  \varepsilon^{\gamma_1} t(\varepsilon) + \xi(\varepsilon),
\end{equation} 
where 
$\xi(\varepsilon)/ \varepsilon^{\gamma_1} t(\varepsilon)  \rightarrow 0$ in $L^1$ as $\varepsilon \downarrow 0$. 
The proof of (\ref{annk}) relies on the following lemma, which provides the asymptotics of the transition times
between different surfaces.
\begin{lemma} \label{lennn}
For each $1 \leq k \leq m$, there is a constant $c_k > 0$ such that 
\[
\lim_{\varepsilon \downarrow 0} (\varepsilon^{\gamma_k} \mathrm{E} \sigma^{x, \varepsilon}_1 )  = c_k,
\]
uniformly in $x \in S_{l,k}$, $l \in \mathbb{Z}^d$. Moreover, the random variables
$\sigma^{x, \varepsilon}_1/ \mathrm{E} \sigma^{x, \varepsilon}_1$ are uniformly integrable in  $\varepsilon > 0$ and $x \in S_{l,k}$, $l \in \mathbb{Z}^d$. 
\end{lemma}
This lemma can be proved by considering the times of excursions from $S_{l,k}$ to the set $\Gamma_\varkappa(S_{l,k})$ (defined
as in Section~\ref{neigh})  and from $\Gamma_\varkappa(S_{l,k})$ to $\partial D$. The former were studied in Lemma~\ref{mltime},
while the latter are much shorter when $\varepsilon$ is small. The number of such excursions prior to  $\sigma^{x, \varepsilon}_1$ depends
on $\varkappa$, and its asymptotics, as $\varkappa \downarrow 0$, can be derived using the arguments similar to those in Section~\ref{neigh}.
The uniform integrability claimed here requires the uniform integrability statements from   Lemmas~\ref{lee2} and 
Lemma~\ref{mltime}. 

Formula (\ref{annk}) follows from Lemma~\ref{lennn} since the number of renewal events for 
$({X}^{ x, \varepsilon}_{{\sigma}^{{x, \varepsilon}}_{n}},  
{\sigma}^{{x, \varepsilon}}_{n})$ prior to time $t$ grows nearly linearly with $t$, with the coefficient equal to 
$(\int_{M} \mathrm{E} \sigma^{x, \varepsilon}_1 d \mu^\varepsilon(x))^{-1}$, where $\mu^\varepsilon$ is the invariant measure 
of the chain $Z^{x,\varepsilon}_n$ on $M$. 

From (\ref{lln1b}), (\ref{cltdb}), (\ref{annk}), and the proximity of $Y^{x,\varepsilon}_{N(\varepsilon)}$ and $X^{x, \varepsilon}_{t(\varepsilon)}$, it follows that
\[
\lim_{\varepsilon \downarrow 0}  \frac{ \mathrm{E} X^{x,\varepsilon}_{t(\varepsilon)}}{c  \varepsilon^{\gamma_1}  
t(\varepsilon)}  = 
\tilde{a};~~~ \lim_{\varepsilon \downarrow 0}   
\frac{   X^{x,\varepsilon}_{t(\varepsilon)}}{c  \varepsilon^{\gamma_1}  
t(\varepsilon)}  = \tilde{a}~~{\rm in}~{\rm probability},
\]
\[
\frac{X^{x,\varepsilon}_{t(\varepsilon)} -  \mathrm{E} X^{x,\varepsilon}_{t(\varepsilon)} }{\sqrt{ c  \varepsilon^{\gamma_1}  
t(\varepsilon)}} 
\rightarrow N(0, \tilde{B} )~~~{\rm in}~~{\rm distribution},~~
{\rm as}~~ {\varepsilon \downarrow 0}. 
\]
The result now follows with $a  = c \tilde{a}  $ and $B = c \tilde{B}$. \qed
\\
\\

Now, let us  briefly discuss the case when $\gamma_1 < 0$, i.e., all the surfaces are repelling. 
Assuming that the process $X^{x,\varepsilon}_t$ starts at $x \in D$, it behaves similarly to $X^{x,\varepsilon}_t$ at time scales
$t(\varepsilon) \ll \varepsilon^{\gamma_1}$. The process $X^{x,\varepsilon}_t$  satisfies the Law of Large Numbers with some drift $a \in \mathbb{R}^d$ and the Central Limit Theorem with a non-degenerate
diffusion matrix $B$. Therefore, we should expect similar behavior for $X^{x,\varepsilon}_t$. However, if $x \in D_{l,k}$ for some $l \in
\mathbb{R}^d$, $1 \leq k \leq m$, then $X^{x,\varepsilon}_t$ will not escape from $D_{l,k}$ prior to time $t(\varepsilon)$ (it takes
time of order $\varepsilon^{\gamma_k}$ to reach $S_{l,k}$, as follows from Lemmas \ref{prre} and \ref{lee1} (part (c))). Thus $X^{x,\varepsilon}_t$ process will remain bounded with probability close to one in this case. 

At longer time scales, i.e., $\varepsilon^{\gamma_k}  \ll t(\varepsilon) \ll \varepsilon^{\gamma_{k+1}}$ for some $1 \leq k < m$,
the process $X^{x,\varepsilon}_t$ behaves similarly to $X^{x}_t$ while the former remains in $D$. However, $X^{x,\varepsilon}_t$
can make excursions into the domains $x \in D_{l,k'}$ with $k' \leq k$ (it takes time of order $\varepsilon^{\gamma_{k'}}$ to
reach $D_{l,k'}$ if the process starts at $x \in D$ and time of the same order to leave a small neighborhood of $D_{l,k'}$ and return to
the diffusive behavior in $D$). Thus, at such time scales,  $X^{x,\varepsilon}_t$,  with $x \in D$, will behave diffusively, but will be slowed down, compared
to $X^{x,\varepsilon}_t$, by a constant factor, due to a positive proportion of time spent inside the domains $D_{l,k'}$ with 
$l \in \mathbb{R}^d$, $1 \leq k' \leq k$. 

Finally, at the longest time scales, $t(\varepsilon) \gg \varepsilon^{\gamma_{m}}$, the process $X^{x,\varepsilon}_t$ will reach $D$ in
time that is negligible compared to $t(\varepsilon)$ for every $x \in \mathbb{R}^d$. Thus the behavior will be diffusive for every initial
point. We thus have the following result, which we provide without a formal proof.
\begin{theorem} Suppose that $\gamma_1  < 0$. There    is  a vector $a \in \mathbb{R}^d$, a positive-definite symmetric 
matrix $B$ and positive constants $1 = c_0 \geq c_1 \geq ... \geq c_m$ such that, for $\varepsilon^{\gamma_k}  \ll t(\varepsilon) \ll \varepsilon^{\gamma_{k+1}}$ (where $\gamma_0: = 0$ and $\gamma_{m+1}: = -\infty$) and $x \in D$ (or $x \in \mathbb{R}^d$ if $k = m$),
\[
\lim_{\varepsilon \downarrow 0}  \frac{ \mathrm{E} X^{x,\varepsilon}_{t(\varepsilon)}}{ 
t(\varepsilon)}  = c_k a;~~~ \lim_{\varepsilon \downarrow 0} \frac{   X^{x,\varepsilon}_{t(\varepsilon)}}{ 
t(\varepsilon)}  = c_k a~~~in~probability,
\]
\[
\frac{X^{x,\varepsilon}_{t(\varepsilon)} -\mathrm{E} X^{x,\varepsilon}_{t(\varepsilon)}}{\sqrt{ t(\varepsilon)  }} \rightarrow 
N(0, c_kB)~~~in~~distribution,~~
as~~\varepsilon \downarrow 0. 
\]
\end{theorem}
\noindent
\\
\\
\noindent {\bf \large Acknowledgments}: The work of  L. Koralov was supported by the Simons Foundation Fellowship (award number 678928).
\\
\\


\begin{thebibliography}{999999}

\bibitem{Ba} Y. Bakhtin, {\it Noisy heteroclinic networks}, Probability Theory and Related Fields, 150 (1): 1--42,   2011.

\bibitem{BCP} Y. Bakhtin, H-B. Chen, Z. Pajor-Gyulai, {\it Rare Transitions in Noisy Heteroclinic Networks}, preprint.

\bibitem{BL}  V. Betz, S. Le Roux, {\it Multi-scale metastable dynamics and the asymptotic stationary distribution of
perturbed Markov chains}, Stoch. Process. Appl. 126 (11), 3499--3526 (2016).


\bibitem{D1} M. Day, {\it Recent progress on the small parameter exit problem},
Stochastics 20 (1987), no. 2, 121--150.

\bibitem{D2} M. Day, {\it Boundary local time and small parameter exit problems with characteristic boundaries},
SIAM J. Math. Analysis, 20 (1989), no. 1, 222--248.

\bibitem{D3} M. Day, {\it Mathematical approaches to the problem of noise-induced exit}, Stochastic
analysis, control, optimization and applications, 269--287, Systems Control
Found. Appl., Birkhauser, Boston, MA, 1999.
\bibitem{Fich} G. Fichera, {\it On a unified theory of boundary value
problems for elliptic-parabolic equations of second order},
Matematika, 1963, Volume 7, Issue 6, 99--122.


\bibitem{F85}  M.I. Freidlin, {\it Functional Integration and Partial Differential Equations}, Princeton
University Press, 1985.

\bibitem{Finf} M. I. Freidlin, {\it Long time influence of small perturbations and motion on the simplex of invariant probability measures}, to appear in
Pure and Applied Functional Analysis.


\bibitem{FK-MC} M. I. Freidlin, L. Koralov, {\it Metastable distributions of Markov chains with rare transitions}. 
J. Stat. Phys. 167 (2017), no. 6, 1355–-1375.

\bibitem{FK21} M. I. Freidlin, L. Koralov, {\it Asymptotics in the Dirichlet Problem 
for Second Order Elliptic Equations with Degeneration on the Boundary}, Journal of Differential Equations,  332 (2022), 202–218.

\bibitem{FW} M. I. Freidlin, A. D. Wentzell, {\it Random
Perturbations of Dynamical Systems}, Springer 2012.

\bibitem{Has} R. Z. Hasminskii, {\it Diffusion processes and elliptic equations degenerating at the boundary of a region} (Russian),
 Teor. Veroyatnost. i Primenen., 3,  1958, pp 430--451.

\bibitem{H2} R. Z. Hasminskii, {\it The averaging principle for parabolic and elliptic differential equations and Markov processes with small diffusion} (Russian), Teor. Verojatnost. i Primenen., 8, 1963, 3--25. 

\bibitem{Kato} T. Kato, {\it Perturbation Theory for Linear Operators}, 2nd ed., 1995, Springer.

\bibitem{LX} C. Landim, T. Xu, {\it Metastability of finite state Markov chains: a recursive procedure
to identify slow variables for model reduction}, ALEA Lat. Am. J. Probab.
Math. Stat. 13 (2016), no. 2, 725--751.





%





















%


\bibitem{MS} B. J. Matkowsky, Z. Schuss, {\it The exit problem for randomly perturbed dynamical
systems}, SIAM J. Appl. Math., 33 (1977), pp. 365--382.




\bibitem{Pin} R. G.  Pinsky, {\it Positive Harmonic Functions and Diffusion}, Cambridge University Press, 1995.





\bibitem{OR} E. V. Radkevich, {\it Equations with nonnegative characteristics form. I.}, J. Math. Sci. (N.Y.) 158 (2009), no. 3, 297–452.


\bibitem{OR2} E. V. Radkevich, {\it Equations with nonnegative characteristic form. II.}, J. Math. Sci. (N.Y.) 158 (2009), no. 4, 453–604.

\bibitem{Va} B. R. Vainberg, V. V. Grushin, {\it Uniformly nonelliptic problems. II.}, (Russian)
Mat. Sbornik 73 (115), 1967, 126–154.


\end{thebibliography}
\end{document}